\pdfoutput=1
\RequirePackage{ifpdf}
\ifpdf 
\documentclass[pdftex]{sigma}
\else
\documentclass{sigma}
\fi

\usepackage{dsfont}
\usepackage{float}
\usepackage[sc,small]{caption}
\usepackage{bm}
\usepackage{enumerate}
\usepackage{tikz}
\usetikzlibrary{arrows, backgrounds}
\usepackage{bbm}
\usepackage[mathscr]{eucal}
\usepackage{rotating}

\newtheorem{thm}{Theorem}
\newtheorem{conv}[thm]{Convention}

\newtheorem{lem}[thm]{Lemma}
\newtheorem{prop}[thm]{Proposition}
\theoremstyle{definition}
\newtheorem{Example}[thm]{Example}
\newtheorem{defi}[thm]{Definition}

\newtheorem{rem}[thm]{Remark}

\newcommand\derive[1]{\ell^{\scriptscriptstyle{\lceil} #1 \scriptscriptstyle{\rceil}}}
\newcommand\e[2]{e_{\scriptscriptstyle{#1}}^{\scriptscriptstyle{#2}}}
\newcommand\er[2]{p_{\scriptscriptstyle{#1}}^{\scriptscriptstyle{#2}}}
\newcommand\nilpoh[1]{h_{\scriptscriptstyle{\lfloor} #1 \scriptscriptstyle{\rfloor}}}
\newcommand\nilpot[1]{\ell_{\scriptscriptstyle{\lfloor} #1 \scriptscriptstyle{\rfloor}}}
\newcommand\cc[2]{c_{#1}^{\scriptscriptstyle(#2 \scriptscriptstyle)}}
\newcommand\tracen[2]{{\tr}_{{#1}} (#2)}
\newcommand\traceN[2]{{\tr}_{{#1}}^{} (#2)}
\newcommand\Tracen[2]{{\tr}_{{#1}} \big(#2 \big)}
\renewcommand{\r}[2]{r_{\scriptscriptstyle{#1}}^{\scriptscriptstyle{#2}}}

\DeclareMathOperator{\End}{End}
\DeclareMathOperator{\Hom}{Hom}
\DeclareMathOperator{\id}{id}
\DeclareMathOperator{\tr}{tr}
\DeclareMathOperator{\Img}{Im}

\newcommand{\retract}[3]{
\draw[fill] (#1,#2) -- (#1+#3,#2) -- (#1+ #3/2,#2+#3) -- cycle;}
\newcommand{\embed}[3]{
\draw[fill] (#1,#2) -- (#1+#3,#2) -- (#1+ #3/2,#2-#3) -- cycle;}
\newcommand{\sovereign}[2]{
\draw[fill=white] (#1-0.1,#2-0.1) rectangle (#1+0.1,#2+0.1);}

\usetikzlibrary{decorations.markings}
\tikzset{
  ma/.style={
    decoration={markings,mark=at position 0.99999 with {\arrow[scale=1.65,#1]{>}}},
    postaction={decorate},
    shorten >=0.4pt}}	
\tikzset{
  ba/.style={
    decoration={markings,mark=at position 0.99999 with {\arrow[scale=1.85,#1]{>}}},
    postaction={decorate},
    shorten >=0.4pt}}
\tikzset{
  iba/.style={
    decoration={markings,mark=at position 0.00001 with {\arrow[scale=1.85,#1]{<}}},
    postaction={decorate},
    shorten >=0.4pt}}	    	
\tikzset{
  bam/.style={
    decoration={markings,mark=at position 0.5 with {\arrow[scale=1.85,#1]{>}}},
    postaction={decorate},
    shorten >=0.4pt}}
\tikzset{
  ibam/.style={
    decoration={markings,mark=at position 0.5 with {\arrow[scale=1.85,#1]{<}}},
    postaction={decorate},
    shorten >=0.4pt}}

\begin{document}

\numberwithin{equation}{section}
\numberwithin{thm}{section}

\allowdisplaybreaks

\renewcommand{\thefootnote}{$\star$}

\renewcommand{\PaperNumber}{017}

\FirstPageHeading

\ShortArticleName{On the Killing form of Lie Algebras in Symmetric Ribbon Categories}

\ArticleName{On the Killing form of Lie Algebras\\
in Symmetric Ribbon Categories\footnote{This paper is a~contribution to the Special Issue on New Directions in Lie
Theory.
The full collection is available at
\href{http://www.emis.de/journals/SIGMA/LieTheory2014.html}{http://www.emis.de/journals/SIGMA/LieTheory2014.html}}}

\Author{Igor BUCHBERGER and J\"urgen FUCHS}
\AuthorNameForHeading{I.~Buchberger and J.~Fuchs}

\Address{Teoretisk fysik, Karlstads Universitet, Universitetsgatan 21, S--65188 Karlstad, Sweden}

\Email{\href{mailto:igor.buchberger@kau.se}{igor.buchberger@kau.se}, \href{mailto:juerfuch@kau.se}{juerfuch@kau.se}}
\URLaddress{\url{http://www.ingvet.kau.se/juerfuch/}}

\ArticleDates{Received September 30, 2014, in f\/inal form February 20, 2015; Published online February 26, 2015}

\Abstract{As a~step towards the structure theory of Lie algebras in symmetric monoidal categories we establish results
involving the Killing form.
The proper categorical setting for discussing these issues are symmetric ribbon categories.}


\Classification{17Bxx; 18D35; 18D10; 18E05}

\renewcommand{\thefootnote}{\arabic{footnote}}
\setcounter{footnote}{0}

\section{Introduction}

At times, a~mathematical notion reveals its full nature only after viewing it in a~more general context than the one in
which it had originally appeared, and often this is achieved by formulating it in the appropriate categorical framework.
To mention a~few examples, the octonions are really an instance of an \emph{associative} commutative algebra, namely
once they are regarded as an object in the category of ${\ensuremath{\mathbb Z}}_2^{\times3}$-graded vector spaces with
suitably twisted associator~\cite{alMaj,baez16}; chiral algebras in the sense of Beilinson and Drinfeld can be viewed as
Lie algebras in some category of $\mathcal D$-modules~\cite{soib4}; and vertex algebras (which are Beilinson--Drinfeld
chiral algebras on the formal disk) are singular commutative associative algebras in a~certain functor
category~\cite{borC22}.
As another illustration, not only do Hopf algebras furnish a~vast generalization of group algebras, but indeed a~group
\emph{is} a~Hopf algebra, namely a~Hopf algebra in the category of sets.

In this note we study aspects of Lie algebras from a~categorical viewpoint.
This already has a~long tradition, see e.g.~\cite{bafm,cofW,fiMo, gure,gure4}.
For instance, it is well known that Lie superalgebras and, more generally, Lie color algebras, are ordinary Lie algebras
in suitable categories, see Example~\ref{ex:lialgs} below.
Our focus here is on the structure theory of Lie algebras in symmetric monoidal categories which, to the best of our
knowledge, has not been investigated in purely categorical terms before.
We restrict our attention to a~few particular issues, the main goal being to identify for each of them an adequate
categorical setting that both allows one to develop the respective aspect closely parallel with the classical situation
and at the same time covers an as broad as possible class of cases.
In particular we do not assume the underlying category to be semisimple, or Abelian, or linear (i.e., enriched over
vector spaces), or that the tensor product functor is exact, albeit large classes of examples that are prominent in the
literature share (some of) these properties.

On the other hand we insist on imposing relevant conditions directly on the underlying category.
They are thereby not necessarily the weakest possible, and in fact they could often be considerably weakened to
conditions on subclasses of objects.
For instance, even an object in a~non-braided category can be endowed with a~Lie algebra structure, provided that it
comes with a~self-invertible Yang--Baxter operator in the sense of~\cite[Definition~2.6]{goVe} that takes over the role
of the self-braiding.
But this way one potentially loses the clear distinction between categorical and non-categorical aspects.
Indeed, a~benef\/it of the abstraction inherent in the categorical point of view is that it allows one to neatly separate
features that apply only to a~subclass of examples from those which are essential for the concepts and results in
question and are thereby generic.
We are specif\/ically interested in the proper notion of Lie subalgebra and ideal, and in the relevance of the Killing form.
Recall that non-degeneracy of the Killing form furnishes a~criterion for semisimplicity of a~Lie algebra over the
complex numbers, but ceases to do so for e.g.\ Lie algebras over a~f\/ield of non-zero characteristic or for Lie
superalgebras.

The rest of this note is organized as follows.
Section~\ref{sec2} is devoted to basic aspects of Lie algebras in monoidal categories.
In Section~\ref{sec2_1} we present the def\/inition of a~Lie algebra in a~symmetric additive category and illustrate it by
examples.
We then introduce the notions of nilpotent and solvable Lie algebras in Section~\ref{sec2_2}.
In Section~\ref{sec2_3} we brief\/ly describe subalgebras and ideals of Lie algebras, paying attention to the dif\/ferences
that result from the specif\/ication of the relevant class of monics.
In Section~\ref{sec3} we study aspects of the Killing form.
To this end we f\/irst review in Section~\ref{sec3_1} some properties of the partial trace.
In Section~\ref{sec3_2} we give the def\/inition of the Killing form of a~Lie algebra in a~symmetric ribbon category and
show that it is symmetric, while invariance is established in Section~\ref{sec3_3}.
Finally in Section~\ref{sec3_4} we explore the role of non-degeneracy of the Killing form.

\section{Lie algebras}
\label{sec2}

We freely use pertinent concepts from category theory, and in particular for algebras in monoidal categories.
Some of these are recalled in the Appendix.
The product of a~Lie algebra is antisymmetric and satisf\/ies the Jacobi identity.
This requires that morphisms can be added and that there is a~notion of exchanging the factors in a~tensor product.
Moreover, it turns out that one can proceed in full analogy to the classical case only if the latter is
a~\emph{symmetric braiding}.
Thus we make the

\begin{conv}
\label{conv:addsymm}
In the sequel, unless stated otherwise, by a~category we mean an additive symmetric monoidal category, with symmetric
braiding~$c$.
\end{conv}

To our knowledge, Lie algebras in this setting were f\/irst considered in~\cite{gure,gure4}.
To present the def\/inition, we introduce the short-hand notation
\begin{gather}
\cc U n:= c^{}_{U_{}^{\otimes(n-1)}, U} \in {\ensuremath{{\End}_{\ensuremath{\mathcal C}}}}\big(U^{\otimes n}\big)
\label{ccUn}
\end{gather}
for multiple self-braidings, as well as $\nilpoh{2}:= h$ and
\begin{gather}
\nilpoh{n}:= h \circ (\id_U \otimes \nilpoh{n-1}) \in {\ensuremath{{\Hom}_{\ensuremath{\mathcal C}}}}\big(U^{\otimes n},U\big)
\qquad
\text{for}
\quad
n > 2
\label{nilpoh}
\end{gather}
for iterations of a~morphism $h \in {\ensuremath{{\Hom}_{\ensuremath{\mathcal C}}}}(U\otimes U,U)$.

\subsection{Lie algebras in additive symmetric monoidal categories}
\label{sec2_1}

With the notations~\eqref{ccUn} and~\eqref{nilpoh} we have

\begin{defi}
A~\emph{Lie algebra} in a~category {\ensuremath{\mathcal C}}\ is a~pair $(L, \ell)$ consisting of an object $L\in
{\ensuremath{\mathcal C}}$ and a~morphism $\ell_{} \in {\ensuremath{{\Hom}_{\ensuremath{\mathcal C}}}}(L\otimes L,L)$
that satisf\/ies
\begin{gather}
\label{antisym}
\nilpot 2 \circ \big(\id_L^{ \otimes 2} + \cc L 2\big) = 0
\end{gather}
(antisymmetry) and
\begin{gather}
\label{jacobi}
\nilpot 3 \circ \big[\id_L^{ \otimes 3} + \cc L 3 +{(\cc L 3)}^2\big] = 0
\end{gather}
(Jacobi identity).
\end{defi}

Here and below we assume monoidal categories to be strict.
The additional occurrences of the associativity constraint in the non-strict case are easily restored (the explicit
expression can e.g.\ be found in~\cite[Def\/inition~2.1]{goVe}). By abuse of terminology, also the object~$L$ is called
a~Lie algebra; the morphism~$\ell$ is referred to as the \emph{Lie bracket} of~$L$.

\begin{Example}\label{ex:lialgs}\quad
\begin{enumerate}\itemsep=0pt
\item[(i)] For $\ensuremath{\Bbbk}$ a~f\/ield and ${\ensuremath{\mathcal C}} = \ensuremath{\mathscr V{\rm ect}}_{\ensuremath{\Bbbk}}$
the category of \ensuremath{\Bbbk}-vector spaces, with the symmetric braiding given by the f\/lip $v\otimes v' {\mapsto}
v' \otimes v$, we recover ordinary Lie algebras.

\item[(ii)] A~Lie superalgebra over $\ensuremath{\Bbbk}$~\cite{frKa,kac7} is a~Lie algebra in the category $\mathscr
S\ensuremath{\mathscr V{\rm ect}}_{\ensuremath{\Bbbk}}$ of ${\ensuremath{\mathbb Z}}_2$-graded \ensuremath{\Bbbk}-vec\-tor
spaces with the braiding given by the superf\/lip, acting as $v\otimes v' {\mapsto} (-1)^{\deg(v) \deg(v')}(v'\otimes v)$
on homogeneous elements (and extended by bilinearity).

\item[(iii)] A~Lie color algebra (or color Lie algebra)~\cite{BMpz,riWy2, scheu2} is a~Lie algebra in the
category~$\Gamma$-$\ensuremath{\mathscr V{\rm ect}}$ of~$\Gamma$-graded vector spaces, where~$\Gamma$ is a~f\/inite Abelian
group endowed with a~skew bicha\-rac\-ter~$\varphi$.
The braiding acts on homogeneous elements analogously as for superalgebras, i.e.\ as a~f\/lip multiplied by a~phase factor
$\varphi(\deg(v),\deg(v'))$.

For $\Gamma = {\ensuremath{\mathbb Z}}_2$ with the unique cohomologically nontrivial skew bicharacter this yields Lie
superalgebras, while for~$\varphi$ cohomologically trivial it yields ordinary Lie algebras for any~$\Gamma$.

\item[(iv)] A~so-called Hom-Lie algebra~$L$ in a~symmetric monoidal category {\ensuremath{\mathcal C}}, for which the
Jacobi identity is deformed with the help of an automorphism of~$L$ (f\/irst introduced in~\cite{haLSi} for the case
${\ensuremath{\mathcal C}} = \ensuremath{\mathscr V{\rm ect}}_{\ensuremath{\Bbbk}}$) is a~Lie algebra in a~category
$\mathscr H{\ensuremath{\mathcal C}}$ whose objects are pairs consisting of an object~$U$ of {\ensuremath{\mathcal C}}\
and an automorphism of~$U$~\cite{caGo}.
\end{enumerate}
\end{Example}

\begin{rem}\label{rem:brLie}\quad
\begin{enumerate}\itemsep=0pt
\item[(i)] When trying to directly generalize the notion of a~Lie algebra to monoidal categories with a~generic braiding, two
independent Jacobi identities must be considered~\cite{zhzh}.
\item[(ii)] A~more conceptual approach to the case of generic braiding has led to the notion of a~quantum Lie bracket
on an object~$L$, satisfying a~generalized Jacobi identity, which is related to the adjoint action for Hopf algebras in
{\ensuremath{\mathcal C}}\ that are compatible with an additional coalgebra structure on~$L$~\cite{maji34}.
In another approach one deals with a~whole family of~$n$-ary products on an object of~$\Gamma$-\ensuremath{\mathscr
V{\rm ect}}, with~$\Gamma$ a~f\/inite Abelian group endowed with a~(not necessarily skew) bicharacter~$\psi$~\cite{pare25}.
(If~$\psi$ is skew, this reduces to Lie color algebras.) This approach has the advantage of keeping important aspects of
the classical case, e.g.\ the primitive elements of a~Hopf algebra in~$\Gamma$-\ensuremath{\mathscr V{\rm ect}}\ and the
derivations of an associative algebra in~$\Gamma$-\ensuremath{\mathscr V{\rm ect}}\ carry such a~generalized Lie algebra
structure.
There is also a~further generalization~\cite{pare27} to the case that {\ensuremath{\mathcal C}}\ is the category of
Yetter--Drinfeld modules over a~Hopf algebra with bijective antipode.
\end{enumerate}
\end{rem}

\begin{rem}
Given an associative algebra $(A,m)$ in {\ensuremath{\mathcal C}}, the morphism
\begin{gather}
\ell_m:= m - m\circ c_{A,A}
\label{eq:commutator}
\end{gather}
clearly satisf\/ies the def\/ining relations~\eqref{antisym} and~\eqref{jacobi} of a~Lie bracket.
Thus the pair $(A,\ell_m)$ is a~Lie algebra in {\ensuremath{\mathcal C}}; we refer to such a~Lie algebra as
a~\emph{commutator Lie algebra} and to its Lie bracket $\ell_m$ as the commutator of the associative multiplication~$m$.
\end{rem}

\begin{rem}
If the category {\ensuremath{\mathcal C}}\ is in addition Abelian, then one can def\/ine the \emph{universal enveloping
algebra} of a~Lie algebra $L \in {\ensuremath{\mathcal C}}$ as the quotient of the tensor algebra $T(L)$, as described
in Example~\ref{exa:4}(i), by the two-sided ideal $I(L)$ generated by the image of the morphism $\ell - \id_{L\otimes L}
+ c_{L,L}$ (regarded as an endomorphism of $T(L)$).
In many cases of interest, such as for Lie color algebras~\cite{gure4,khar9} or if {\ensuremath{\mathcal C}}\ is the
category of comodules over a~coquasitriangular bialgebra in $\ensuremath{\mathscr V{\rm ect}}$~\cite{fiMo}, the ideal $I(L)$
inherits a~Hopf algebra structure from $T(L)$, so that the quotient is naturally a~Hopf algebra as well.
This still applies~\cite[Section~6]{pare27} in the case of the Lie algebras in categories of Yetter--Drinfeld modules
mentioned in Remark~\ref{rem:brLie}(ii).
\end{rem}

\begin{rem}
The class of all Lie algebras in {\ensuremath{\mathcal C}}\ forms the objects of a~category, with a~morphism~$f$ from
a~Lie algebra $(L,\ell)$ to a~Lie algebra $(L',\ell')$ being a~morphism $f \in {\ensuremath{{\Hom}_{\ensuremath{\mathcal C}}}}(L,L')$ satisfying
\begin{gather}
\label{mirphLie}
\ell' \circ (f\otimes f) = f \circ \ell.
\end{gather}
The category of associative algebras in {\ensuremath{\mathcal C}}\ is def\/ined analogously; supplementing the mapping
from an associative algebra to its commutator Lie algebra by the identity on the sets of algebra morphisms def\/ines
a~functor from the category of associative algebras in {\ensuremath{\mathcal C}}\ to the one of Lie algebras in
{\ensuremath{\mathcal C}}.
Note that, as the condition~\eqref{mirphLie} is nonlinear, these categories are not additive.
\end{rem}

\subsection{Nilpotent and solvable Lie algebras}
\label{sec2_2}

Basic aspects of the structure theory of Lie algebras can be developed in the same way as in the classical case.
We start with

\begin{defi}
An \emph{Abelian} Lie algebra is a~Lie algebra with vanishing Lie product, $\ell = 0$.
\end{defi}

\begin{rem}\quad
\begin{enumerate}\itemsep=0pt
\item[(i)] Trivially, $(U,\ell=0)$ is an Abelian Lie algebra for any object $U \in {\ensuremath{\mathcal C}}$.
A~generic object cannot be endowed with any other Lie algebra structure.
\item[(ii)] The commutator Lie algebra $(A,\ell_m)$ built from a~commutative associative algebra $(A,m)$ is Abelian.
\end{enumerate}
\end{rem}

Next we introduce the concepts of nilpotent and solvable Lie algebras.
Recall that according to the def\/inition~\eqref{nilpoh}, for a~Lie algebra $(L,\ell)$ the $(n-1)$-fold iteration of the
Lie bracket is denoted by $\nilpot{n} \in {\ensuremath{{\Hom}_{\ensuremath{\mathcal C}}}}(L^{\otimes n}, L)$.
We now introduce in addition morphisms $\derive{n} \in {\ensuremath{{\Hom}_{\ensuremath{\mathcal C}}}}(L^{\otimes
2^{(n-1)}}, L)$~by
\begin{gather*}
\derive{2}:= \ell
\qquad
\text{and}
\qquad
\derive{n}:=\ell\circ \big(\derive{n-1}\otimes\derive{n-1}\big)
\qquad
\text{for} \ \ n  >  2,
\end{gather*}
respectively.
Then we can give

\begin{defi}
Let $L = (L,\ell)$ be a~Lie algebra in {\ensuremath{\mathcal C}}.
\begin{enumerate}\itemsep=0pt
\item[(i)] $L$ is called \emph{nilpotent} if\/f $\nilpot{n} = 0$ for suf\/f\/iciently large~$n$.
\item[(ii)] $L$ is called \emph{solvable} if\/f $\derive{n} = 0$ for suf\/f\/iciently large~$n$.
\item[(iii)] $L$ is called \emph{derived nilpotent} if\/f $\nilpot{n} \circ \ell^{\otimes n} = 0$ for suf\/f\/iciently
large~$n$.
\end{enumerate}
\end{defi}

It follows directly from the def\/initions that a~nilpotent Lie algebra is solvable.
Indeed,
\begin{gather*}
\derive{n} = \nilpot{n}\circ \big(\derive{n-1}\otimes\derive{n-2}\otimes \dots \otimes \ell \otimes \id_{L}\otimes \id_{L}\big).
\end{gather*}
But we even have

\begin{prop}
A~derived nilpotent Lie algebra is solvable.
\end{prop}

\begin{proof}
We have directly $\derive{3} = \nilpot{2} \circ \ell^{\otimes 2}$, while for $n {>} 3$ one shows by induction that
\begin{gather*}
\derive{n+1} = \nilpot{n} \circ \ell^{\otimes n}\circ \big(\derive{n-1}\otimes\derive{n-1} \otimes
\derive{n-2}\otimes\derive{n-2} \otimes \dots \otimes \ell \otimes \ell \otimes \id_{L}^{ \otimes 4}\big).
\end{gather*}
Thus $\nilpot{n} \circ \ell^{\otimes n} = 0$ implies $\derive{n+1} = 0$, so that the claim follows directly from the
def\/initions of derived nilpotency and of solvability.
\end{proof}

\subsection{Subalgebras and ideals}
\label{sec2_3}

A subalgebra of an algebra~$A$ is an isomorphism class of monics to~$A$ that are algebra morphisms.
For spelling out this notion in detail, one must specify the class of monics one is considering.
We present the cases that either all monics are admitted, corresponding to working with subobjects, or only split monics
are admitted, corresponding to retracts.
In the latter case we give

\begin{defi}\label{subalgebra}\quad
\begin{enumerate}\itemsep=0pt
\item[(i)] A~\emph{retract subalgebra}, or \emph{subalgebra retract}, of an algebra $(A,m)$ in {\ensuremath{\mathcal C}}\ is
a~retract $(B,\e{B}{A},\r{A}{B})$ of~$A$ for which
\begin{gather}
m \circ \big(\e{B}{A}\otimes \e{B}{A}\big) = \er{A}{B} \circ m \circ \big(\e{B}{A}\otimes \e{B}{A}\big),
\label{eq:subalg-m}
\end{gather}
where $\er AB = \e{B}{A} \circ \r{A}{B}$ is the idempotent associated with the retract.
\item[(ii)] A~\emph{retract Lie subalgebra}, or \emph{Lie subalgebra retract}, of a~Lie algebra $(L,\ell)$ in
{\ensuremath{\mathcal C}}\ is a~retract $(K,\e{K}{L},\r{L}{K})$ of~$L$ for which
\begin{gather}
\ell \circ \big(\e{K}{L}\otimes \e{K}{L}\big) = \er{L}{K} \circ \ell \circ \big(\e{K}{L}\otimes \e{K}{L}\big).
\label{eq:subalg-l}
\end{gather}
\end{enumerate}
\end{defi}

Note that part (ii) of Def\/inition~\ref{subalgebra} is redundant~-- it is merely a~special case of part~(i), as we do not
specify any properties of the product~$m$ in~(i).
For clarity we present the Lie algebra version nevertheless separately; it is worth recalling that in the Lie algebra
case the Convention~\ref{conv:addsymm} is in ef\/fect.
Also note that owing to $p \circ e = e$,~\eqref{eq:subalg-m} is equivalent to
\begin{gather*}
m \circ \big(\er{A}{B}\otimes \er{A}{B}\big) = \er{A}{B} \circ m \circ \big(\er{A}{B}\otimes \er{A}{B}\big),
\end{gather*}
and analogously for~\eqref{eq:subalg-l}, as well as for~\eqref{eq:lieideal} below.

\begin{defi}
\label{ideal}
A~\emph{retract ideal}, or \emph{ideal retract}, of a~Lie algebra $(L,\ell)$ in {\ensuremath{\mathcal C}}\ is a~retract
$(K,\e{K}{L},\r{L}{K})$ of~$L$ for which
\begin{gather}
\ell \circ \big(\e{K}{L}\otimes \id_{L}\big) = \er{L}{K} \circ \ell \circ \big(\e{K}{L}\otimes \id_{L}\big).
\label{eq:lieideal}
\end{gather}
\end{defi}

We call a~Lie algebra~$L$ indecomposable if\/f it is non-Abelian and its only retract ideals are the zero object and~$L$
itself.
Put dif\/ferently, we have

\begin{defi}
\label{sindecomposable}
An \emph{indecomposable} Lie algebra is a~Lie algebra which is not Abelian and which does not possess any non-trivial retract ideal.
\end{defi}

By antisymmetry of~$\ell$, the equality~\eqref{eq:lieideal} is equivalent to $\ell \circ (\id_{L} \otimes \e{K}{L}) =
\er{L}{K} \circ \ell \circ (\id_{L} \otimes \e{K}{L})$.
For an associative algebra, the equality analogous to~\eqref{eq:lieideal} instead def\/ines a~\emph{left} ideal, while
exchanging the roles of the two morphisms $\id$ and~$e$ gives the notion of right ideal.

If $(K,e,r)$ is a~non-trivial retract of a~Lie algebra $(L, \ell)$, then together with $p = r \circ e$ also $p':= \id_{L}-p$
is a~non-zero idempotent, and hence if {\ensuremath{\mathcal C}}\ is idempotent complete, then there is
a~non-trivial retract $(K',e',r')$ of~$L$ such that $p' = r' \circ e'$ and $L \cong  K \oplus K'$ as an object in
{\ensuremath{\mathcal C}}.
In the present context a~particularly interesting case is that both~$K$ and $K'$ are retract ideals of~$L$.
To account for this situation we give

\begin{defi}
Given a~f\/inite family of Lie algebras $(K_i, \ell_i)$, the \emph{direct sum} Lie algebra of the family is the Lie
algebra $(L, \ell)$ whose underlying object~$L$ is the direct sum $L:= \bigoplus_i K_i$ of the objects underlying the
Lie algebras $(K_i, \ell_i)$, and whose Lie bracket is given~by
\begin{gather}
\label{ds_bracket}
\ell:= \sum\limits_i \e{K_i}{L} \circ \ell_i \circ \big(\r{L}{K_i}\otimes \r{L}{K_i}\big).
\end{gather}
\end{defi}

Note that since $\r{L}{K_i}\circ \e{K_j}{L} = \delta_{i,j} \id_{K_j}$, for a~direct sum Lie algebra we have $\ell \circ
(\er{L}{K_i} \otimes \er{L}{K_j})=0$ for $i {\ne} j$.
Moreover, it follows that $\ell\circ(\er{L}{K_i} \otimes \id_{L}) = \er{L}{K_i}\circ\ell\circ(\er{L}{K_i} \otimes \id_{L})$
for every~$i$, i.e.\ the Lie algebras $(K_i, \ell_i)$ are retract ideals of the direct sum Lie algebra~$L$.
Conversely, one has

\begin{lem}
If a~Lie algebra $(L, \ell)$ is, as an object, the direct sum $L = \bigoplus_i K_i$, and each of the retracts $K_i$ is
a~retract ideal of~$L$, then $(L, \ell)$ is the direct sum of the $K_i$ also as a~Lie algebra.
\end{lem}

\begin{proof}
Because of $\id_{L} = \sum\limits_i \er{K_i}{L}$ the Lie bracket~$\ell$ can be written as the triple sum
\begin{gather}
\label{ijk}
\ell = \sum\limits_{i,j,k} \er{L}{K_k} \circ \ell\circ \big(\er{L}{K_i} \otimes \er{L}{K_j}\big).
\end{gather}
Now since $K_i$ and $K_j$ are ideals, one has $\er{L}{K_k} \circ \ell\circ (\er{L}{K_i} \otimes \er{L}{K_j}) =
\delta_{i,k} \delta_{j,k} \er{L}{K_k} \circ \ell\circ (\er{L}{K_i} \otimes \er{L}{K_j})$.
As a~consequence the expression~\eqref{ijk} reduces to
\begin{gather*}
\ell = \sum\limits_{i} \e{K_i}{L} \circ \ell_i \circ \big(\r{L}{K_i} \otimes \r{L}{K_i}\big)
\qquad
\text{with}
\qquad
\ell_i:= \r{L}{K_i} \circ \ell \circ \big(\e{K_i}{L} \otimes \e{K_i}{L}\big).
\end{gather*}
Hence there are brackets $\ell_i$ satisfying~\eqref{ds_bracket}.
Moreover, invoking the uniqueness of direct sum decomposition in the additive category {\ensuremath{\mathcal C}}, it
follows (notwithstanding the non-additivity of the category of Lie algebras in {\ensuremath{\mathcal C}}) that these
brackets coincide with the assumed Lie brackets on the retract ideals $K_i$ up to isomorphism of Lie algebras.
\end{proof}

Next, consider a~retract ideal $(K,e,r)$ of a~Lie algebra~$L$ in an idempotent-complete category.
We call~$K$ a~\emph{primitive} retract ideal if\/f the idempotent $p = e \circ r$ cannot be written as a~sum $p = p_1 +
p_2$ of two non-zero idempotents such that at least one of the corresponding retracts $(K_1,e_1,r_1)$ and
$(K_2,e_2,r_2)$ is a~retract ideal of~$L$.
Now if the retracts associated with two idempotents in ${\ensuremath{{\End}_{\ensuremath{\mathcal C}}}}(L)$, say~$p$ and
$p'$, are retract ideals, then so is the retract associated with the idempotent $p\circ p' \in
{\ensuremath{{\End}_{\ensuremath{\mathcal C}}}}(L)$.
If~$p$ is in addition primitive, then either $p\circ p' = p$ or $p\circ p' = 0$.
A~primitive retract ideal is hence an indecomposable Lie algebra.

It is worth pointing out that in the Def\/initions~\ref{subalgebra} and~\ref{ideal} there is no need to assume the
existence of an image of the Lie bracket~$\ell$.
In contrast, if we spell out the corresponding def\/initions for the case of subobjects, we do have to make such an
assumption:

\begin{defi}\quad
Let $(L,\ell)\in {\ensuremath{\mathcal C}}$ be a~Lie algebra for which the image $\Img(\ell)$ exists in
{\ensuremath{\mathcal C}}.
\begin{enumerate}\itemsep=0pt
\item[(i)] A~\emph{subobject Lie subalgebra}, or \emph{Lie subalgebra subobject}, of~$L$ is a~subobject $(K,\e{K}{L})$
of~$L$ for which the image $ \Img\big(\ell \circ (\e{K}{L}\otimes \e{K}{L}) \big)$ is a~subobject of~$K$.
\item[(ii)] A~\emph{subobject ideal}, or \emph{ideal subobject}, of~$L$ is a~subobject $(K,\e{K}{L})$ of~$L$ for which
the image $ \Img\big(\ell \circ (\e{K}{L}\otimes \id_{L}) \big)$ is a~subobject of~$K$.
\end{enumerate}
\end{defi}

Accordingly, Def\/inition~\ref{sindecomposable} then gets replaced~by

\begin{defi}
A~\emph{simple} Lie algebra is a~Lie algebra which is not Abelian and which does not possess any non-trivial subobject
ideal.
\end{defi}

\begin{Example}\label{ex:..M}\qquad
\begin{enumerate}\itemsep=0pt
\item[(i)] Given a~unital associative algebra $(A,m,\eta)$ in {\ensuremath{\mathcal C}}, the commutator
$\ell_m$~\eqref{eq:commutator} satisf\/ies $\ell_m \circ (\eta \otimes \id_A) \equiv  m \circ (\eta \otimes \id_A) - m \circ
c_{A,A} \circ (\eta \otimes \id_A) = \id_A - \id_A = 0$.
Thus, provided that $({\bf1},\eta)$ is a~subobject of~$A$~-- meaning that the unit morphism~$\eta$ is monic, which is
true under rather weak conditions~\cite{macdJ}, e.g.\ if {\ensuremath{\mathcal C}}\ is \ensuremath{\Bbbk}-linear with
\ensuremath{\Bbbk}\ a~f\/ield~\cite[Section~7.7]{deli}~-- $({\bf1},\eta)$ is an Abelian subobject ideal of the commutator Lie
algebra $(A,\ell_m)$.
\item[(ii)] If there exists in addition a~morphism $\epsilon\in {\ensuremath{{\Hom}_{\ensuremath{\mathcal
C}}}}(A,{\bf1})$ such that $\xi:= \epsilon \circ \eta \in {\ensuremath{{\End}_{\ensuremath{\mathcal C}}}}({\bf1})$ is
invertible, then $({\bf1},\eta,\tilde\epsilon)$ with $\tilde\epsilon:= \xi^{-1} \circ \eta$ is a~retract of~$A$, and thus
even an Abelian retract ideal of the Lie algebra $(A,\ell_m)$, and we have a~direct sum decomposition
\begin{gather}
\label{AA'}
A= {\bf1} \oplus A'
\end{gather}
of Lie algebras.
\end{enumerate}
\end{Example}

\begin{Example}\label{ex:M..}\quad
\begin{enumerate}\itemsep=0pt
\item[(i)] If~$A$ carries the structure of a~Frobenius algebra in {\ensuremath{\mathcal C}}, then its counit~$\varepsilon$ is
a~natural candidate for the morphism~$\epsilon$ that is assumed in Example~\ref{ex:..M}(ii).
If, moreover, the Frobenius algebra~$A$ is strongly separable in the sense of Def\/inition~\ref{def:Frob}(ii) as well as
symmetric, then $\varepsilon \circ \eta = {\ensuremath{{\dim}_{\ensuremath{\mathcal C}}}}(A)$ equals the dimension
of~$A$, i.e.\ $({\bf1},\eta,\tilde\varepsilon)$ being a~retract ideal of the commutator Lie algebra $(A,\ell_m)$ is
equivalent to~$A$ having invertible dimension.
\item[(ii)] As a~special case, consider the tensor product $A_U:= U \otimes U^\vee$ of an object~$U$ with its right
dual in a~(strictly) sovereign category {\ensuremath{\mathcal C}}.
The object $A_U$ carries a~canonical structure of a~symmetric Frobenius algebra
$(A_U,m_U,\eta_U,\Delta_U,\varepsilon_U)$, for which all structural morphisms are simple combinations of the left and
right evaluation and coevaluation morphisms for~$U$, namely
\begin{gather}
\label{AU}
m_U = \id_U\otimes d_U\otimes\id_{U^\vee_{}},
\qquad
\eta_U = b_U
\end{gather}
and similarly for the coproduct $\Delta_U$ and counit $\varepsilon_U$.
We refer to $A_U$ with this algebra and co\-al\-gebra structure as a~(full) \emph{matrix algebra}.
Any matrix algebra\ is Morita equivalent to the trivial Frobenius algebra ${\bf1}$, and we have $\varepsilon_U \circ
\eta_U = {\ensuremath{{\dim}_{\ensuremath{\mathcal C}}}}(U)$, so ${\bf1}$ is a~retract ideal of the commutator Lie
algebra $(A_U,\ell_{m_U})$ if\/f~$U$ has invertible dimension.
\item[(iii)] A~simple realization of the situation described in~(ii) is obtained by taking {\ensuremath{\mathcal C}}\
to be the category of f\/inite-dimensional complex representations of a~classical f\/inite-dimensional simple Lie algebra
$\mathfrak L$ in $\ensuremath{\mathscr V{\rm ect}}_{{\mathbbm C}}$ (or in $\ensuremath{\mathscr V{\rm ect}}_\ensuremath{\Bbbk}$,
with \ensuremath{\Bbbk}\ any algebraically closed f\/ield of characteristic zero) and $U = D \in {\ensuremath{\mathcal
C}}$ to be the def\/ining representation of $\mathfrak L$.
The retract $A_D'$ in the resulting direct sum decomposition $A_D = {\bf1}  \oplus  A_D'$ is then a~simple Lie algebra
if\/f $\mathfrak L = \mathfrak{sl}_n({\mathbbm C})$, in which case $A_D' = \mathrm{A d}$ is (isomorphic to) the adjoint
representation of $\mathfrak{sl}_n({\mathbbm C})$ and $A_D = {\bf1}  \oplus  \mathrm{A d}$ is just the familiar
decomposition $\mathfrak{gl}_n({\mathbbm C}) \cong  {\mathbbm C}  \oplus  \mathfrak{sl}_n({\mathbbm C})$.
\item[(iv)] Another situation in which the retract $A_U'$ is a~simple Lie algebra is the case that
{\ensuremath{\mathcal C}}\ is the category of f\/inite-dimensional complex representations of either the Mathieu group
$\mathbb M = \mathbb M_{23}$ or of $\mathbb M = \mathbb M_{24}$, and~$U$ is the 45-dimensional irreducible
representation of~$\mathbb M$; $A_D'$ is then the 2024-dimensional irreducible $\mathbb M$-representation.
(The decomposition $A_D = {\bf1} \oplus  A_D'$ is easily checked by hand from the character table of $\mathbb M$, which
can e.g.\ be found in Tables~V and~VI of~\cite{jameG}.)
\end{enumerate}
\end{Example}

\begin{Example}
Let $B = (B,m,\eta,\Delta,\varepsilon)$ be a~bialgebra and $(B,\ell = m - m\circ c)$ its commutator Lie algebra.
Consider a~subobject $(P,e)$ of~$B$ satisfying
\begin{gather}
\label{P,e}
\Delta \circ e = e \otimes \eta + \eta \otimes e.
\end{gather}
By the bialgebra axiom it follows directly that $\Delta \circ m \circ (e \otimes e) = (m \circ (e \otimes e)) \otimes \eta
+ 2 e \otimes e + \eta \otimes (m \circ (e \otimes e))$, and similarly for $\Delta \circ m \circ c \circ (e \otimes e)$.
Hence the morphism $\ell_P:= \ell \circ(e \otimes e)$ obeys $ \Delta \circ \ell_P = \ell_P \otimes \eta + \eta \otimes
\ell_P$.
In particular, if $(P,e)$ is maximal with the property~\eqref{P,e}, then it is a~subobject Lie subalgebra of $(B,\ell)$.
\end{Example}

\begin{rem}
In applications, the relevance of a~Lie algebra arises often through its linear representations.
In the present note we do not dwell on representation theoretic issues, except for the following simple observation.
For $A_U = U \otimes U^\vee$ the associative algebra with product~$m_U$ as in~\eqref{AU}, the object~$U$ is naturally an
$A_U$-module, with representation morphism $\rho:= \id_U \otimes d_U \in {\ensuremath{{\Hom}_{\ensuremath{\mathcal
C}}}}(A_U\otimes U,U)$.
The representation property is verif\/ied as follows:
\begin{gather*}
\rho_U \circ (\id_{A_U} \otimes \rho_U) \equiv (\id_U \otimes d_U) \circ (\id_{A_U} \otimes \id_U \otimes d_U) =
\id_U \otimes d_U \otimes d_U
\\
\phantom{\rho_U \circ (\id_{A_U} \otimes \rho_U)}{}
= (\id_U \otimes d_U) \circ (\id_U \otimes d_U \otimes \id_{U^\vee_{}} \otimes \id_U) \equiv \rho_U \circ (m_U\otimes \id_U).
\end{gather*}
Being a~module over the associative algebra $(A_U,m_U)$, $U$ is also a~module over the commutator Lie algebra
$(A_U,\ell_{m_U})$ and thus, in case~$U$ has invertible dimension, also over the Lie algebra~$A_U'$ appearing in the
direct sum decomposition~\eqref{AA'}.
It then follows in particular that if every object of {\ensuremath{\mathcal C}}\ is a~retract of a~suitable tensor power
of~$U$, then each object of {\ensuremath{\mathcal C}}\ carries a~natural structure of an $A_U'$-module, whereby
{\ensuremath{\mathcal C}}\ becomes equivalent to a~subcategory of the category $A_U'\mbox{-mod}_{\ensuremath{\mathcal C}}$ of $A_U'$-modules in {\ensuremath{\mathcal C}}.
In the special case of $U = D$ being the def\/ining representation of $\mathfrak{sl}_n({\mathbbm C})$,
{\ensuremath{\mathcal C}}\ is in fact equivalent to $A_U'\mbox{-mod}_{\ensuremath{\mathcal C}}$ as a~monoidal category.
\end{rem}

\section{The Killing form}\label{sec3}

\subsection{Partial traces}\label{sec3_1}

The Killing form of a~Lie algebra in $\ensuremath{\mathscr V{\rm ect}}_\ensuremath{\Bbbk}$ is def\/ined as a~trace, actually as a~partial trace.
Canonical traces, and thus partial traces, of morphisms exist in a~category {\ensuremath{\mathcal C}}\ if\/f
{\ensuremath{\mathcal C}}\ is a~full monoidal subcategory of a~ribbon category~\cite{josv}.
Accordingly from now on we restrict our attention to ribbon categories (but for now do not yet impose our
Convention~\ref{conv:addsymm}).

One may think of the basic data of a~(strictly monoidal) ribbon category to be the right and left dualities, $b$, $d$ and
$\tilde b$, $\tilde d$, the braiding~$c$, and the sovereign structure~$\sigma$, rather than the twist~$\theta$.
Then instead of expressing the sovereign structure with the help of the twist as in~\eqref{defsigma}, conversely the
twist is expressed as
\begin{gather}
\label{thetavssigma}
\theta_U = (\id_U \otimes \tilde d_U) \circ (c_{U,U} \otimes \sigma_U) \circ (\id_U \otimes b_U).
\end{gather}
It will prove to be convenient to exploit the graphical calculus for (strict) monoidal categories and for algebras
therein that has been developed repeatedly in dif\/ferent guises, see e.g.~\cite{coDu2,joSt5,KAss,maji34,seli3,yama12} or
Appendix~A~of~\cite{fjfrs}.
For instance, in this graphical description the def\/ining property~\eqref{fvvf} of the sovereign structure takes the form
\begin{gather*}
\begin{tikzpicture}
[scale=0.65] \node [above] at (-1.5,1.5) {$\scriptstyle{{}^{\vee}U}$}; \draw (-1.5, 0.6)-- (-1.5,1.6); \draw (-2, 0)
rectangle (-1, 0.6); \node at (-1.5, 0.3) {$\scriptstyle{{}^{\vee }f}$}; \draw (-1.5, -1.6)-- (-1.5,0); \node [below] at
(-1.5,-1.5) {$\scriptstyle{V^{\vee}}$}; \node [right] at (-1.5,-0.75) {$\scriptstyle{\sigma_{\scriptscriptstyle{V}}}$};
\sovereign{-1.5}{-0.75}  \node at (0.25, 0){$=$};  \draw (2, 0.2)-- (2,1.6); \draw (1.5, -0.4) rectangle (2.5, 0.2);
\node [above] at (2,1.5) {$\scriptstyle{{}^{\vee}U}$}; \node at (2.1, -0.1) {$\scriptstyle{f^{\vee}}$}; \sovereign
{2}{0.95} \node [right] at (2,0.95) {$\scriptstyle{\sigma_{\scriptscriptstyle{U}}}$}; \draw (2, -0.4)-- (2,-1.6); \node
[below] at (2,-1.5) {$\scriptstyle{V^{\vee}}$};
\end{tikzpicture}
\end{gather*}
while the expression~\eqref{thetavssigma} for the twist becomes
\begin{gather*}
\begin{tikzpicture}[scale=0.5]
 \node at (-2.5,0) {$\theta_U$};
 \node at (-1.5,0) {$=$};
 \draw(-0.5,-1.5) to [out=45,in=275 ] (0,0);
 \draw[ma] (0,0) arc [radius=0.5, start angle=180, end angle= 90];
 \draw (0.5,0.5) arc [radius=0.5, start angle=90, end angle= 0];
 \draw[ma] (1,0) arc [radius=0.5, start angle=0, end angle= -90];
 \draw (0.5,-0.5) arc [radius=0.5, start angle=-90, end angle= -130];
 \draw(-0.1, -0.2) to [out=130,in=275 ] (-0.5, 1.5);
 \sovereign{1}{0}
 \node [below] at (-0.5,-1.5) {$\scriptstyle{U}$};
 \node [right] at (1,0) {$\scriptstyle{\sigma_{\scriptscriptstyle{U}}}$};
 \node [above] at (-0.5,1.5) {$\scriptstyle{U}$};
 \node at (2.75,0) {$.$};
\end{tikzpicture}
\end{gather*}

The (left and right) \emph{trace} of an endomorphism $f\in {\ensuremath{{\End}_{\ensuremath{\mathcal C}}}}(V)$ in
a~ribbon category is given by $\tracen{}f = \tilde{d}_V \circ (f \otimes \sigma_{V}) \circ b_{V} = {d}_V \circ
(\sigma_V^{-1} \otimes f) \circ \tilde{b}_{V}
\in {\ensuremath{{\End}_{\ensuremath{\mathcal C}}}}({\bf1})$; pictorially,
\begin{gather*}
\begin{tikzpicture}[scale=0.6]
 \node at (-2.8,-0.75) {$\tr(f) $};
 \node at (-1.4,-0.75) {$=$};
 \draw (0, -0.45)-- (0,0);
 \draw (-0.5, -1.05) rectangle (0.5,-0.45) ;
 \draw (0, -1.05)-- (0,-1.5);
 \node at (0,-0.75) {$\scriptstyle{f} $};
 \draw[ba] (0,0) arc [radius=0.75, start angle=180, end angle= 90];
 \draw (0.75,0.75) arc [radius=0.75, start angle=90, end angle= 0];
 \draw (1.5, 0)-- (1.5,-1.5);
 \draw[ba] (1.5,-1.5) arc [radius=0.75, start angle=0, end angle= -90];
 \draw (0.75,-2.25) arc [radius=0.75, start angle=-90, end angle= -180];
 \node [right] at (1.5,-0.75) {$\scriptstyle{\sigma_{\scriptscriptstyle{V}}}$};
 \sovereign{1.5}{-0.75};
 \node at (3.1,-0.75) {$=$};
 \draw (5, 0)-- (5,-1.5);
 \draw[ba] (5,-1.5) arc [radius=0.75, start angle=-180, end angle= -90];
 \draw (5.75,-2.25) arc [radius=0.75, start angle=-90, end angle= 0];
 \node [right] at (3.75,-0.65) {$\scriptstyle{\sigma^{-1}_{\scriptscriptstyle{V}}}$};
 \draw (6.5, -0.45)-- (6.5,0);
 \draw (6, -1.05) rectangle (7,-0.45) ;
 \draw (6.5, -1.05)-- (6.5,-1.5);
 \node at (6.5,-0.75) {$\scriptstyle{f} $};
 \draw[ba] (6.5,0) arc [radius=0.75, start angle=0, end angle= 90];
 \draw (5.75,0.75) arc [radius=0.75, start angle=90, end angle= 180];
 \sovereign{5}{-0.75};
 \node at (7.75,-0.75) {$.$};
\end{tikzpicture}
\end{gather*}
For def\/ining the Killing form we also need a~generalization, the notion of a~\emph{partial trace}: for a~morphism $f \in
{\ensuremath{{\Hom}_{\ensuremath{\mathcal C}}}}(U_1 {\otimes}\cdots{\otimes} U_n \otimes V, V)$ we introduce the partial
trace with respect to~$V$~by
\begin{gather*}
\tracen{n+1}f:= \tilde{d}_V \circ (f \otimes \sigma_V) \circ (\id_{U_1\otimes\dots\otimes U_n}^{} {\otimes} b_{V}) \in
{\ensuremath{{\Hom}_{\ensuremath{\mathcal C}}}}(U_1 {\otimes}\cdots{\otimes} U_n,{\bf1}).
\end{gather*}
For $n = 1$ this looks graphically as follows:
\begin{gather}
\label{tracen2}
\begin{split}
& \begin{tikzpicture}[scale=0.6]
 \node at (-3.4,-0.75) {$\tracen{2}f $};
 \node at (-1.9,-0.75) {$=$};
 \draw (0, -0.45)-- (0,0);
 \draw (-1, -1.05) rectangle (0.5,-0.45) ;
 \node at (-0.25,-0.75) {$\scriptstyle{f} $};
 \draw (-0.5, -1.05)-- (-0.5,-2.25);
 \draw (0, -1.05)-- (0,-1.5);
 \node[below] at (-0.5, -2.25) {$\scriptstyle{U}$};
 \draw[ba] (0,0) arc [radius=0.75, start angle=180, end angle= 90];
 \draw (0.75,0.75) arc [radius=0.75, start angle=90, end angle= 0];
 \draw (1.5, 0)-- (1.5,-1.5);
 \draw[ba] (1.5,-1.5) arc [radius=0.75, start angle=0, end angle= -90];
 \draw (0.75,-2.25) arc [radius=0.75, start angle=-90, end angle= -180];
 \node [right] at (1.5,-0.75) {$\scriptstyle{\sigma_{\scriptscriptstyle{V}}}$};
 \node at (5.2,-0.75) {$\in~{\ensuremath{{\Hom}_{\ensuremath{\mathcal C}}}}(U,{\bf1})$};
 \sovereign{1.5}{-0.75};
\end{tikzpicture}
\end{split}
\end{gather}
for $f \in {\ensuremath{{\Hom}_{\ensuremath{\mathcal C}}}}(U\otimes V,V)$.
A~fundamental property of the partial trace~\eqref{tracen2} is the following:

\begin{lem}
For any morphism $f \in {\ensuremath{{\Hom}_{\ensuremath{\mathcal C}}}}(U\otimes V,V)$ in a~ribbon category one has
\begin{gather}
\tracen{2}{f} \circ \theta_U = \tracen{2}f.
\label{eqtheta}
\end{gather}
\end{lem}

\begin{proof}
This follows directly by invoking functoriality of the twist:
\begin{gather*}
\tracen{2}{f} \circ \theta_U = \Tracen{2}{\theta_V^{} \circ f \circ \big(\id_U \otimes \theta_V^{-1}\big)} = \Tracen{2}{f \circ
\big[\id_U \otimes \big(\theta_V^{-1} \circ \theta_V^{}\big)\big]} = \tracen{2}f.
\end{gather*}
Alternatively one may `drag the morphism~$f$ along the~$V$-loop'; graphically:
\begin{gather}
\label{prooftheta}
\begin{split}
& \begin{tikzpicture}[scale=0.53]
 \node at (-3.4,-0.75) {$\tracen{2}f $};
 \node at (-1.9,-0.75) {$\equiv$};
 \draw (0, -0.45)-- (0,0);
 \draw (-1, -1.05) rectangle (0.5,-0.45) ;
 \node at (-0.25,-0.75) {$\scriptstyle{f} $};
 \draw (-0.5, -1.05)-- (-0.5,-2.5);
 \draw (0, -1.05)-- (0,-1.5);
 \node[below] at (-0.5, -2.5) {$\scriptstyle{U}$};
 \draw[ba] (0,0) arc [radius=0.75, start angle=180, end angle= 90];
 \draw (0.75,0.75) arc [radius=0.75, start angle=90, end angle= 0];
 \draw (1.5, 0)-- (1.5,-1.5);
 \draw[ba] (1.5,-1.5) arc [radius=0.75, start angle=0, end angle= -90];
 \draw (0.75,-2.25) arc [radius=0.75, start angle=-90, end angle= -180];
 \sovereign{1.5}{-0.75};
 \node at (2.5,-0.75) {$=$};
 \draw (4, -1.5)-- (4,0);
 \draw (5, -1.05) rectangle (6.5,-0.45) ;
 \node at (5.75,-0.75) {$\scriptstyle{{}^{\vee\!}\!f} $};
 \draw[ba] (4,0) arc [radius=0.75, start angle=180, end angle= 90];
 \draw (4.75,0.75) arc [radius=0.75, start angle=90, end angle= 0];
 \draw (5.5, 0)-- (5.5,-0.45);
 \draw (5.5, -1.05)-- (5.5,-1.5);
 \draw[ba] (5.5,-1.5) arc [radius=0.75, start angle=0, end angle= -90];
 \draw (4.75,-2.25) arc [radius=0.75, start angle=-90, end angle= -180];
 \draw (6, 0.5)-- (6,-0.45);
 \draw[ba](3.5,0.5) arc [radius=1.25, start angle=180, end angle= 90];
 \draw(4.75,1.75) arc [radius=1.25, start angle=90, end angle= 0];
 \draw (3.5, 0.5)-- (3.5,-2.5);
 \sovereign{5.5}{-1.45};
 \node at (7.5,-0.75) {$=$};
 \draw (9, -1.5)-- (9,0);
 \draw (10, -1.05) rectangle (11.5,-0.45) ;
 \node at (10.75,-0.75) {$\scriptstyle{f}^{\vee}$};
 \draw[ba] (9,0) arc [radius=0.75, start angle=180, end angle= 90];
 \draw (9.75,0.75) arc [radius=0.75, start angle=90, end angle= 0];
 \draw (10.5, 0)-- (10.5,-0.45);
 \draw (10.5, -1.05)-- (10.5,-1.5);
 \draw[ba] (10.5,-1.5) arc [radius=0.75, start angle=0, end angle= -90];
 \draw (9.75,-2.25) arc [radius=0.75, start angle=-90, end angle= -180];
 \draw (11, 0.5)-- (11,-0.45);
 \draw[ba] (8.5,0.5) arc [radius=1.25, start angle=180, end angle= 90];
 \draw (9.75,1.75) arc [radius=1.25, start angle=90, end angle= 0];
 \draw (8.5, 0.5)-- (8.5,-2.5);
 \sovereign{10.5}{0};
 \sovereign{11}{0};
 \node at (12.5,-0.75) {$=$};
 \draw (-0.25+15, -0.45)-- (-0.25+15,0);
 \draw (-1+15, -1.05) rectangle (0.5+15,-0.45) ;
 \node at (-0.25+15,-0.75) {$\scriptstyle{f} $};
 \draw (-0.5+15, -1.05)-- (-0.5+15,-1.55);
 \draw (0+15, -1.05)-- (0+15,-1.5);
 \draw[ba] (-0.25+15,0) arc [radius=0.75+0.125, start angle=180, end angle= 90];
 \draw (0.75-0.125+15,0.75+0.125) arc [radius=0.75+0.125, start angle=90, end angle= 0];
 \draw (1.5+15, 0)-- (1.5+15,-1.5);
 \draw[ba] (1.5+15,-1.5) arc [radius=0.75, start angle=0, end angle= -90];
 \draw (0.75+15,-2.25) arc [radius=0.75, start angle=-90, end angle= -180];
 \draw[ba] (2+15,-1.5) arc [radius=1.25, start angle=0, end angle= -90];
 \draw (0.75+15,-2.75) arc [radius=1.25, start angle=-90, end angle= -180];
 \draw (2+15, -1.5) -- (2+15,0) ;
 \draw(2+15,0.) arc [radius=1.5+0.125, start angle=0, end angle= 90];
 \draw [iba] (0.5-0.125+15,1.5+0.125) arc [radius=1.5+0.125, start angle=90, end angle= 180];
 \draw (2-3.25+15, -2.95) -- (2-3.25+15,0) ;
 \sovereign{1.5+15}{-0.75};
 \sovereign{2+15}{-0.75};
 \node at (3+15,-0.75) {$=$};
 \draw (4.75+15, -0.45)-- (4.75+15,0);
 \draw (4+15, -1.05) rectangle (5.5+15,-0.45) ;
 \node at (4.75+15,-0.75) {$\scriptstyle{f} $};
 \draw (5+15, -1.05)-- (5+15,-1.5);
 \draw[ba] (4.75+15,0) arc [radius=0.75+0.125, start angle=180, end angle= 90];
 \draw (5.75-0.125+15,0.75+0.125) arc [radius=0.75+0.125, start angle=90, end angle= 0];
 \draw (6.5+15, 0)-- (6.5+15,-1.5);
 \draw[ba] (6.5+15,-1.5) arc [radius=0.75, start angle=0, end angle= -90];
 \draw (5.75+15,-2.25) arc [radius=0.75, start angle=-90, end angle= -180];
 \sovereign{6.5+15}{-0.75};
 \draw (4.5+15, -1.05) to [out=270,in=90] (4.5+15-0.45,-1.9);
 \draw[->] (5.0+15-0.35,-2.) arc [radius=0.25, start angle=0, end angle= -90];
 \draw (4.75+15-0.35,-2.25) arc [radius=0.25, start angle=-90, end angle= -155];
 \draw (5+15-0.35,-2) arc [radius=0.25, start angle=0, end angle= 90];
 \draw [->](4.4+15-0.35,-2.25) to [out=90,in=180] (4.75+15-0.35, -1.75);
 \draw (4.4+15-0.35,-2.25) -- (4.4+15-0.35,-2.95) ;
 \sovereign{5+15-0.35}{-2};
\end{tikzpicture}
\end{split}\!\!\!\!\!\!\!
\end{gather}
Here the f\/irst and third equalities hold by def\/inition of the left and right dual morphisms ${}^{\vee\!\!}f$ and $f^\vee$,
respectively (and by monoidality of~$\sigma$), the second equality is sovereignty, and the f\/inal equality follows by the
functoriality of the braiding.
\end{proof}

By the same sequence of manipulations one obtains analogous identities for other partial traces.
In particular we have the following generalization of the cyclicity property of the ordinary trace:

\begin{lem}
\label{lem:cyclic}
Partial traces are cyclic up to braidings and twists: For any pair of morphisms $f \in
{\ensuremath{{\Hom}_{\ensuremath{\mathcal C}}}}(U\otimes X,Y)$ and $g \in {\ensuremath{{\Hom}_{\ensuremath{\mathcal
C}}}}(V\otimes Y,X)$ in a~ribbon category one has
\begin{gather}
\Tracen{3} {f \circ [\id_U \otimes g]} = \Tracen{3} {g \circ [\id_V \otimes f] \circ [(c_{U,V}^{} \circ (\theta_U \otimes
\id_V))\otimes \id_X]}.
\label{eqcyclic}
\end{gather}
\end{lem}

\begin{rem}
Instead of dragging, as in~\eqref{prooftheta}, the morphism~$f$ along the~$V$-loop in clockwise direction, we could
equally well drag it in counter-clockwise direction.
In the situation considered in Lemma~\ref{lem:cyclic}, this amounts to the identity
\begin{gather*}
\Tracen{3} {f \circ [\id_U \otimes g]} = \Tracen{3} {g \circ [\id_V \otimes f] \circ \big[\big(c_{V,U}^{-1} \circ \big(\id_U \otimes
\theta_V^{-1}\big)\big) \otimes \id_X\big]}
\end{gather*}
instead of~\eqref{eqcyclic}.
This is, however, equivalent to~\eqref{eqcyclic}, as is seen by combining it with the result~\eqref{eqtheta}, as applied
to the morphism $f \circ (\id_U \otimes g)$, and using that $\theta_{U\otimes V} = (\theta_U \otimes \theta_V) \circ
c_{U,V} \circ c_{V,U}$.
\end{rem}

\begin{Example}
\label{exa:svect1}
The category \ensuremath{\mathscr{SV}{\rm ect}}\ of super vector spaces over a~f\/ield \ensuremath{\Bbbk}\ is
\ensuremath{\Bbbk}-linear rigid monoidal and is endowed with a~braiding given by the super-f\/lip map.
\ensuremath{\mathscr{SV}{\rm ect}}\ is semisimple with, up to isomorphims, two simple objects, namely the tensor unit ${\bf1}$
which is the {\bf1}-dimensional vector space \ensuremath{\Bbbk}\ in degree zero, and an object ${S}$ given by \ensuremath{\Bbbk}\ in degree one.
There are two structures of ribbon category on this braided rigid monoidal category.
For the f\/irst, the twist obeys (besides $\theta_{\bf1} = \id_{\bf1}$, which is true in any ribbon category) $\theta_{S}
= -\id_{S}$, which after identifying ${}^{\vee } S = S^\vee$ results in a~strictly sovereign structure and in $\dim(S)=1$.
For the second ribbon structure instead the twist is trivial, $\theta_{S} = \id_{S}$, while the sovereign structure is
non-trivial, namely $\sigma_S = -\id_{S^\vee_{}}$ (upon still identifying ${}^{\vee } S = S^\vee$), and the trace is
a~\emph{supertrace}, in particular now $\dim(S) = -1$.
This second ribbon structure arises from the conventions used in the bulk of the literature on associative and Lie
superalgebras (see e.g.\ \cite{kac7}), while in many applications of super vector spaces in supersymmetric quantum f\/ield
theory the f\/irst ribbon structure is implicit.
\end{Example}

\subsection{The Killing form of a~Lie algebra in a~symmetric ribbon category}\label{sec3_2}

We now turn to morphisms involving Lie brackets.
Accordingly in the sequel our Convention~\ref{conv:addsymm} will again be in ef\/fect.
In particular, as the braiding of {\ensuremath{\mathcal C}}\ is now symmetric, the twist of {\ensuremath{\mathcal C}}\
squares to the identity, and $\theta \equiv  \id$ is an allowed twist.

In a~symmetric ribbon category, we give

\begin{defi}
\label{def:Killing}
The \emph{Killing form} of a~Lie algebra $(L,\ell)$ is the morphism $\kappa \in
{\ensuremath{{\Hom}_{\ensuremath{\mathcal C}}}}(L\otimes L,{\bf1})$ given by the partial trace
\begin{gather*}
\kappa:= \traceN{3}{\theta_L\circ\nilpot{3}}.
\end{gather*}
\end{defi}
Writing out the partial trace, this expression reads more explicitly
\begin{gather*}
\kappa = \tilde{d}_{L} \circ \big([\theta_L\circ\nilpot{3}]\otimes \sigma_L \big) \circ \big(\id_{L}^{ \otimes
2} \otimes b_{L} \big).
\end{gather*}
In the sequel we will use the graphical notation
\begin{gather*}
\begin{tikzpicture}[scale=0.5]
 \node at (-2,0.25) {$\ell$};
 \node at (-1,0.25) {$=$};
 \draw(0,-1) -- (0,0);
 \draw (0,0) arc [radius=0.5, start angle=180, end angle= 0];
 \draw(1,-1) -- (1,0);
 \draw(0.5,0.5) -- (0.5,1.5);
 \draw [fill] (0.40,0.40) rectangle (0.60,0.60);
 \node [above] at (0.5,1.5) {$\scriptstyle{L}$};
 \node [below] at (0,-1) {$\scriptstyle{L}$};
 \node [below] at (1,-1) {$\scriptstyle{L}$};
 \node at (7,0.25) {$\theta_L$};
 \node at (8,0.25) {$=$};
 \draw(9,-1) -- (9,1.5);
 \draw [fill] (9, 0.25) circle [radius=0.15] ;
 \node [above] at (9,1.5) {$\scriptstyle{L}$};
 \node [below] at (9,-1) {$\scriptstyle{L}$};
\end{tikzpicture}
\end{gather*}
for the Lie bracket and the twist of~$L$; then the Killing form is depicted as
\begin{gather*}
\begin{tikzpicture}[scale=0.35]
 \node at (-3.5,0) {$\kappa ~~=$};
 \draw(0,-2) -- (0,0);
 \draw (0,0) arc [radius=0.5, start angle=180, end angle= 0];
 \draw(1,-1) -- (1,0);
 \draw(0.5,0.5) -- (0.5,1.5);
 \draw(-0.25,2.25) -- (-0.25,3);
 \draw [fill] (0.40,0.40) rectangle (0.60,0.60);
 \draw [fill] (-0.35,2.15) rectangle (-0.15,2.35);
 \draw [fill] (-0.25, 2.95) circle [radius=0.175] ;
 \draw (0.5,1.5) arc [radius=0.75, start angle=0, end angle= 180];
 \draw(-1.,1.5) -- (-1,-2);
 \draw (1,-1) arc [radius=0.5, start angle=180, end angle= 270];
 \draw[ba] (2,-1) arc [radius=0.5, start angle=360, end angle= 270];
 \draw(2,-1) -- (2,3);
 \draw[ba] (-0.25,3) arc [radius=1.125, start angle=180, end angle= 90];
 \draw (1.125-0.25,3+1.125) arc [radius=1.125, start angle=90, end angle= 0];
 \node [below] at (-1,-2) {$\scriptstyle{L}$};
 \node [below] at (0,-2) {$\scriptstyle{L}$};
 \draw[fill=white](2-0.15, 1.5-0.15) rectangle (2+0.15, 1.5+0.15);
 \node at (3,0){$.$} ;
\end{tikzpicture}
\end{gather*}
The presence of the twist in Def\/inition~\ref{def:Killing}, which might not have been expected, is needed for the Killing
form to play an analogous role in the structure theory as in the classical case of Lie algebras in $\ensuremath{\mathscr
V{\rm ect}}_\ensuremath{\Bbbk}$.
For comparison we will below also brief\/ly comment on the morphism
\begin{gather}
\label{kappa0}
\kappa_0:= \traceN{3}{\nilpot{3}}.
\end{gather}
which of course coincides with~$\kappa$ if, as in the classical case, $\theta_L =\id_L$.

As illustrated by Example~\ref{exa:svect1}, there are important classes of ribbon categories for which the sovereign
structure is not strict.
Nevertheless, for convenience in the sequel we will suppress the sovereign structure or, in other words, take it to be
strict.
We are allowed to do so because it is entirely straightforward to restore the sovereign structure in all relevant
expressions.

\begin{Example}
For the matrix Lie algebra $(A_U,\ell_{m_U})$, introduced in Example~\ref{ex:M..}(ii) as the commutator Lie algebra of
the matrix algebra $A_U = U \otimes U^\vee$, for any $U\in {\ensuremath{\mathcal C}}$, the Killing form evaluates to
\begin{gather*}
\kappa_{A_U^{}} = 2 \dim_{\theta}(U) \big(\tilde{d}_U \circ [\id_U\otimes d_U \otimes \theta_{U^{\vee}}] \big) - 2
\big(\tilde{d}_U \circ [\id_U\otimes \theta_{U^{\vee}}] \big)^{\otimes 2}
\\
\phantom{\kappa_{A_U^{}}}
\begin{tikzpicture}[scale=0.5]
 \node at (-2.9,0.5) {$= 2\dim_{\theta}(U)$};
\draw[ba] (0,0) arc [radius=1.5, start angle=180, end angle= 90];
\draw (1.5,1.5) arc [radius=1.5, start angle=90, end angle= 0];
\node[below] at (-0.2,0.1) {$\scriptstyle{U}$};
\node[below] at (0.75,0.18) {$\scriptstyle{U^\vee}$};
\node[below] at (2.5,0.1) {$\scriptstyle{U}$};
\node[below] at (3.45,0.18) {$\scriptstyle{U^\vee}$};
\draw [fill] (3-0.08, 0.5) circle [radius=0.1] ;
\draw[ba] (2.5,0) arc [radius=1, start angle=0, end angle= 90];
\draw (1.5,1) arc [radius=1, start angle=90, end angle= 180];
\node at (4.2,0.5) {$-$};
\node at (5.4,0.5) {$2$};
\draw[ba] (5.5+0.4,0.5) arc [radius=0.75, start angle=180, end angle= 90];
\draw (6.25+0.4,1.25) arc [radius=0.75, start angle=90, end angle= 0];
\draw (5.5+0.4,0.5)--(5.5+0.4,0);
\node[below] at (-0.2+6.1,0.1) {$\scriptstyle{U}$};
\node[below] at (1.5+5.7,0.18) {$\scriptstyle{U^\vee}$};
\draw [fill] (7+0.4, 0.4) circle [radius=0.1] ;
\draw (7+0.4,0.5)--(7+0.4,0);
\node[below] at (-0.2+6.1+2.3,0.1) {$\scriptstyle{U}$};
\node[below] at (1.5+5.7+2.3,0.18) {$\scriptstyle{U^\vee}$};
\draw[ba] (7.5+0.6,0.5) arc [radius=0.75, start angle=180, end angle= 90];
\draw (8.25+0.6,1.25) arc [radius=0.75, start angle=90, end angle= 0];
\draw [fill] (7.5+1.5+0.6, 0.4) circle [radius=0.1] ;
\draw (7.5+0.6,0.5)--(7.5+0.6,0);
\draw (9+0.6,0.5)--(9+0.6,0);
\node at (10+0.6,0.5){$,$} ;
\end{tikzpicture}
\end{gather*}
where
\begin{gather*}
\dim_{\theta}(U):= \tilde{d}_U \circ (\theta\otimes \id_U)\circ b_U.
\end{gather*}
Further, recall that if the dimension $\dim(U)$ is invertible, then we have a~direct sum decomposition $A_U = {\bf1}
\oplus  A_U'$.
Using the graphical description
\begin{gather*}
\begin{tikzpicture}[scale=0.6]
\node at (2.5,-1.75){$\id_{U\otimes U^{\vee}} = \e{\ne}{} \circ \r{}{\ne} + \e{A_U'}{} \circ \r{}{A_U'} =
 \displaystyle\frac{b_U\circ\tilde{d}_U}{\dim(U)} + \e{A_U'}{} \circ \r{}{A_U'}= \displaystyle\frac{1}{\dim(U)}$};
\draw[bam] (5+7,-1-2) arc [radius=0.65, start angle=180, end angle= 0];
\draw[bam] (6.3+7,1.5-2) arc [radius=0.65, start angle=360, end angle= 180];
\node at (7+7+0.25,-1.75){$+$};
\draw(7.7+7+0.25,-1-2) arc [radius=0.65, start angle=180, end angle= 0];
\draw (9+7+0.25,1.5-2) arc [radius=0.65, start angle=360, end angle= 180];
\draw[red, ultra thick](7.7+0.65+7+0.25,-1+0.65-2)--(7.7+0.65+7+0.25,1.5-0.65-2);
\node[below] at (8.9+7+0.25,-1.3) {$\scriptscriptstyle{A_U'}$};
\node[below] at (7.7+7+0.25,0.6-3.5) {$\scriptstyle{U}$};
\node[below] at (9.2+7+0.25,0.68-3.5) {$\scriptstyle{U^\vee}$};
\node[below] at (5+7,0.6-3.5) {$\scriptstyle{U}$};
\node[below] at (6.5+7,0.68-3.5) {$\scriptstyle{U^\vee}$};
\end{tikzpicture}
\end{gather*}
of the corresponding idempotents, the Killing form on the ideal $A_U'$ reads
\begin{gather*}
\begin{tikzpicture}[scale=0.5]
 \node at (-12.2,0.25) {$\kappa_{A_U'}^{} = 2 \dim_{\theta}(U) \big(\tilde{d}_U \circ [\id_U\otimes d_U
 \otimes \theta_{U^{\vee}}] \circ [ \e{A_U'}{} \otimes \e{A_U'}{} ] \big)=2 \dim_{\theta}(U)$};
\draw[ba] (0,0) arc [radius=1.5, start angle=180, end angle= 90];
\draw (1.5,1.5) arc [radius=1.5, start angle=90, end angle= 0];
\draw [fill] (3-0.07, 0.5) circle [radius=0.1] ;
\draw[ba] (2.25,0) arc [radius=0.75, start angle=0, end angle= 90];
\draw (1.5,0.75) arc [radius=0.75, start angle=90, end angle= 180];
\draw (0,0) arc [radius=0.75/2, start angle=180, end angle= 360];
\draw (2.25,0) arc [radius=0.75/2, start angle=180, end angle= 360];
\draw[red, ultra thick](0.75/2,-0.75/2)--(0.75/2, -1.2);
\draw[red, ultra thick](2.25+0.75/2,-0.75/2)--(2.25+0.75/2, -1.2);
\node[below] at (0.75/2,-1.05) {$\scriptstyle{A_U'}$};
\node[below] at (2.25+0.75/2,-1.05) {$\scriptstyle{A_U'}$};
\node at (3.75,0){$.$} ;
\end{tikzpicture}
\end{gather*}
\end{Example}

\begin{prop}
The Killing form is \emph{symmetric}, i.e.\ satisfies
\begin{gather}
\label{ksym}
\kappa = \kappa \circ c_{_{L,L}}.
\end{gather}
\end{prop}

\begin{proof}
By functoriality of the twist we can rewrite the Killing form as
\begin{gather}
\label{kappa_alt}
\begin{split}
& \begin{tikzpicture}[scale=0.35]
 \node at (-6.5,0) {$\kappa=\traceN{3}{\theta_L\circ \nilpot{3}}\equiv$};
 \draw(0,-2) -- (0,0);
 \draw (0,0) arc [radius=0.5, start angle=180, end angle= 0];
 \draw(1,-1) -- (1,0);
 \draw(0.5,0.5) -- (0.5,1.5);
 \draw(-0.25,2.25) -- (-0.25,3);
 \draw [fill] (0.40,0.40) rectangle (0.60,0.60);
 \draw [fill] (-0.35,2.15) rectangle (-0.15,2.35);
 \draw [fill] (-0.25, 2.95) circle [radius=0.15] ;
 \draw (0.5,1.5) arc [radius=0.75, start angle=0, end angle= 180];
 \draw(-1.,1.5) -- (-1,-2);
 \draw(1,-1) arc [radius=0.5, start angle=180, end angle= 270];
 \draw[ba] (2,-1) arc [radius=0.5, start angle=360, end angle= 270];
 \draw(2,-1) -- (2,3);
 \draw[ba] (-0.25,3) arc [radius=1.125, start angle=180, end angle= 90];
 \draw (1.125-0.25,3+1.125) arc [radius=1.125, start angle=90, end angle= 0];
 \node [below] at (-1,-2) {$\scriptstyle{L}$};
 \node [below] at (0,-2) {$\scriptstyle{L}$};
\node at (3.3,0){$=$} ;
 \draw(0+5.5,-2) -- (0+5.5,0);
 \draw (0+5.5,0) arc [radius=0.5, start angle=180, end angle= 0];
 \draw(1+5.5,-1) -- (1+5.5,0);
 \draw(0.5+5.5,0.5) -- (0.5+5.5,1.5);
 \draw(-0.25+5.5,2.25) -- (-0.25+5.5,3);
 \draw [fill] (0.40+5.5,0.40) rectangle (0.60+5.5,0.60);
 \draw [fill] (-0.35+5.5,2.15) rectangle (-0.15+5.5,2.35);
 \draw [fill] (1+5.5, -0.5) circle [radius=0.15] ;
 \draw (0.5+5.5,1.5) arc [radius=0.75, start angle=0, end angle= 180];
 \draw(-1.+5.5,1.5) -- (-1+5.5,-2);
 \draw (1+5.5,-1) arc [radius=0.5, start angle=180, end angle= 270];
 \draw[ba] (2+5.5,-1) arc [radius=0.5, start angle=360, end angle= 270];
 \draw(2+5.5,-1) -- (2+5.5,3);
 \draw[ba] (-0.25+5.5,3) arc [radius=1.125, start angle=180, end angle= 90];
 \draw (1.125-0.25+5.5,3+1.125) arc [radius=1.125, start angle=90, end angle= 0];
 \node [below] at (-1+5.5,-2) {$\scriptstyle{L}$};
 \node [below] at (0+5.5,-2) {$\scriptstyle{L}$};
\node at (14.5,0){$\equiv\traceN{3}{ \nilpot{3} \circ( \id_L^{\otimes 2}\otimes \theta_L)}.$} ;
\end{tikzpicture}
\end{split}
\end{gather}
Using the cyclicity~\eqref{eqcyclic} of the partial trace we thus get
\begin{gather}
\label{ksym_f}
\begin{split}
& \begin{tikzpicture}[scale=0.35]
 \node at(-2.5,0){$\kappa =$};
 \draw(0,-2) -- (0,0);
 \draw (0,0) arc [radius=0.5, start angle=180, end angle= 0];
 \draw(1,-1) -- (1,0);
 \draw(0.5,0.5) -- (0.5,1.5);
 \draw(-0.25,2.25) -- (-0.25,3);
 \draw [fill] (0.40,0.40) rectangle (0.60,0.60);
 \draw [fill] (-0.35,2.15) rectangle (-0.15,2.35);
 \draw [fill] (1, -0.5) circle [radius=0.15] ;
 \draw (0.5,1.5) arc [radius=0.75, start angle=0, end angle= 180];
 \draw(-1.,1.5) -- (-1,-2);
 \draw (1,-1) arc [radius=0.5, start angle=180, end angle= 270];
 \draw[ba] (2,-1) arc [radius=0.5, start angle=360, end angle= 270];
 \draw(2,-1) -- (2,3);
 \draw[ba] (-0.25,3) arc [radius=1.125, start angle=180, end angle= 90];
 \draw (1.125-0.25,3+1.125) arc [radius=1.125, start angle=90, end angle= 0];
 \node [below] at (-1,-2) {$\scriptstyle{L}$};
 \node [below] at (0,-2) {$\scriptstyle{L}$};
\node at (3.3,0){$=$} ;
 \draw(-1+5.5,-2) to [out=90,in=275] (0+5.5,0);
 \draw [fill] (-0.5+5.5,-1.05) circle [radius=0.15] ;
 \draw (0+5.5,0) arc [radius=0.5, start angle=180, end angle= 0];
 \draw(1+5.5,-1) -- (1+5.5,0);
 \draw(0.5+5.5,0.5) -- (0.5+5.5,1.5);
 \draw(-0.25+5.5,2.25) -- (-0.25+5.5,3);
 \draw [fill] (0.40+5.5,0.40) rectangle (0.60+5.5,0.60);
 \draw [fill] (-0.35+5.5,2.15) rectangle (-0.15+5.5,2.35);
 \draw [fill] (0.5+5.5, 1.35) circle [radius=0.15] ;
 \draw (0.5+5.5,1.5) arc [radius=0.75, start angle=0, end angle= 180];
 \draw(-1.+5.5,1.5) to [out=275,in=90] (0.5+5.5,-2);
 \draw (1+5.5,-1) arc [radius=0.5, start angle=180, end angle= 270];
 \draw[ba] (2+5.5,-1) arc [radius=0.5, start angle=360, end angle= 270];
 \draw(2+5.5,-1) -- (2+5.5,3);
 \draw[ba] (-0.25+5.5,3) arc [radius=1.125, start angle=180, end angle= 90];
 \draw (1.125-0.25+5.5,3+1.125) arc [radius=1.125, start angle=90, end angle= 0];
 \node [below] at (-1+5.5,-2) {$\scriptstyle{L}$};
 \node [below] at (0.5+5.5,-2) {$\scriptstyle{L}$};
\node at (23.35,0.15){$\!\!\!\!\!\!\!\equiv\traceN{3}{
 [\ell\circ(\id_L\otimes \theta_L)]\circ[\id_L\otimes \ell]\circ[(c_{L,L}\circ(\theta_L\otimes \id_L))\otimes\id_L]}.$};
\end{tikzpicture}
\end{split}\!\!\!\!\!
\end{gather}
The claim now follows by noticing that, via functoriality of the twist, the right hand side of~\eqref{ksym_f} equals
$\traceN{3}{\nilpot{3} \circ (\id_L^{\otimes 2}\otimes \theta_L)} \circ c_{L,L} =\kappa \circ c_{L,L} $.
\end{proof}

The morphism $\kappa_0$ def\/ined in~\eqref{kappa0} satisf\/ies $\kappa_0 = \kappa_0 \circ c_{_{L,L}} \circ (\theta_L\otimes
\id_L)$ instead of~\eqref{ksym}.

\begin{Example}
Recall from Example~\ref{exa:svect1} that there are two ribbon structures on the category
\ensuremath{\mathscr{SV}{\rm ect}}.
For both of them the Killing form~$\kappa$ of a~Lie algebra~$L$ in \ensuremath{\mathscr{SV}{\rm ect}}, i.e.\ a~Lie
superalgebra, is supersymmetric, $\kappa(y,x) = (-1)^{|x| |y|}\kappa(x,y)$ for homogeneous elements $x,y\in L$.
For the ribbon structure that has trivial twist, this trivially holds for $\kappa_0$ as well, while for the one with
non-trivial twist (and strict sovereign structure), $\kappa_0$ is instead symmetric, $\kappa_0(y,x) = \kappa_0(x,y)$.
\end{Example}

\subsection{Invariance of the Killing form}
\label{sec3_3}

We are now going to show that the Killing form is invariant.
The proof relies on the use of the Jacobi identity.
In the graphical description, the Jacobi identity (slightly rewritten using that the braiding is symmetric) reads
\begin{gather*}
\begin{tikzpicture}[scale=0.35]
 \draw(6.75,2.25) -- (6.75,4);
 \draw [fill] (6.65,2.15) rectangle (6.85,2.35);
 \draw (7.5,1.5) arc [radius=0.75, start angle=0, end angle= 180];
 \draw(7.5,1.5) -- (7.5,1);
 \draw(6,1.5) -- (6,-1.5);
 \draw [fill] (7.4,0.9) rectangle (7.6,1.1);
 \draw (8.25,0.25) arc [radius=0.75, start angle=0, end angle= 180];
 \draw(6.75,0.25) -- (6.75,-1.5);
 \draw(8.25,0.25) -- (8.25,-1.5);
 \node[below] at (5.75,-1.3) {$\scriptstyle{L}$};
 \node[below] at (7,-1.3) {$\scriptstyle{L}$};
 \node[below] at (8.25,-1.3) {$\scriptstyle{L}$};
 \node[above] at (6.75,4) {$\scriptstyle{L}$};
 \node at (9.6,1) {$+$};
 \draw(11.75,2.25) -- (11.75,4);
 \draw [fill] (11.65,2.15) rectangle (11.85,2.35);
 \draw (12.5,1.5) arc [radius=0.75, start angle=0, end angle= 180];
 \draw(12.5,1.5) -- (12.5,1);
 \draw(11,1.5) to [out=275,in=150] (11.8, 0.6);
 \draw(11.8,0.6) to [out=330,in=90] (14, -1.5);
 \draw [fill] (12.4,0.9) rectangle (12.6,1.1);
 \draw (13.25,0.25) arc [radius=0.75, start angle=0, end angle= 180];
 \draw(11.75,0.25) -- (11.75,-1.5);
 \draw(13.25,0.25) -- (13.25,-1.5);
 \node[below] at (14.25,-1.3) {$\scriptstyle{L}$};
 \node[below] at (11.5,-1.3) {$\scriptstyle{L}$};
 \node[below] at (12.9,-1.3) {$\scriptstyle{L}$};
 \node[above] at (11.75,4) {$\scriptstyle{L}$};
 \node at (14.9,1) {$+$};
 \draw(17.75,2.25) -- (17.75,4);
 \draw [fill] (17.65,2.15) rectangle (17.85,2.35);
 \draw (18.5,1.5) arc [radius=0.75, start angle=0, end angle= 180];
 \draw(18.5,1.5) -- (18.5,1);
 \draw(17,1.5) -- (17,-1.5);
 \draw [fill] (18.4,0.9) rectangle (18.6,1.1);
 \draw (19.25,0.25) arc [radius=0.75, start angle=0, end angle= 180];
 \draw(17.75,0.25) to [out=275,in=90] (18.25, -1.5);
 \draw(19.25,0.25) to [out=275,in=90] (16.25, -1.5);
 \node[below] at (16,-1.3) {$\scriptstyle{L}$};
 \node[below] at (18.5,-1.3) {$\scriptstyle{L}$};
 \node[below] at (17,-1.3) {$\scriptstyle{L}$};
 \node[above] at (17.75,4) {$\scriptstyle{L}$};
 \node at (22.0,1) {$=\hspace{0.2cm}0.$};
\end{tikzpicture}
\end{gather*}
An alternative version, obtained by further rewriting with the help of antisymmetry, is
\begin{gather}
\begin{split}
& \begin{tikzpicture}[scale=0.35]
 \draw(6.75,2.25) -- (6.75,4);
 \draw [fill] (6.65,2.15) rectangle (6.85,2.35);
 \draw (7.5,1.5) arc [radius=0.75, start angle=0, end angle= 180];
 \draw(7.5,1.5) -- (7.5,1);
 \draw(6,1.5) -- (6,-1.5);
 \draw [fill] (7.4,0.9) rectangle (7.6,1.1);
 \draw (8.25,0.25) arc [radius=0.75, start angle=0, end angle= 180];
 \draw(6.75,0.25) -- (6.75,-1.5);
 \draw(8.25,0.25) -- (8.25,-1.5);
 \node at (9.6,1) {$-$};
 \draw(12.75,2.25) -- (12.75,4);
 \draw [fill] (12.65,2.15) rectangle (12.85,2.35);
 \draw (12,1.5) arc [radius=0.75, start angle= 180, end angle= 0];
 \draw(13.5,1.5) -- (13.5,-1.5);
 \draw(12,1.5) -- (12,1);
 \draw [fill] (11.9,0.9) rectangle (12.1,1.1);
 \draw (12.75,0.25) arc [radius=0.75, start angle=0, end angle= 180];
 \draw(11.25,0.25) -- (11.25,-1.5);
 \draw(12.75,0.25) -- (12.75,-1.5);
 \node at (14.9,1) {$-$};
 \draw(17.75,2.25) -- (17.75,4);
 \draw [fill] (17.65,2.15) rectangle (17.85,2.35);
 \draw (18.5,1.5) arc [radius=0.75, start angle=0, end angle= 180];
 \draw(18.5,1.5) -- (18.5,1);
 \draw(17,1.5) to [out=275,in=90] (17.4,-1.5);
 \draw [fill] (18.4,0.9) rectangle (18.6,1.1);
 \draw (19.25,0.25) arc [radius=0.75, start angle=0, end angle= 180];
 \draw(17.75,0.25) to [out=275,in=90] (16.25, -1.5);
 \draw(19.25,0.25) -- (19.25, -1.5);
 \node at (22.0,1) {$=\hspace{0.2cm}0$};
\end{tikzpicture}
\end{split}
\label{jacobi_pic2}
\end{gather}
(here, as well as in the pictures below, we omit the obvious~$L$-labels).

\begin{prop}
The Killing form is invariant, i.e.
\begin{gather}
\label{kinv}
\kappa \circ (\ell \otimes \id_L) = \kappa \circ (\id_{L}\otimes \ell).
\end{gather}
Graphically,
\begin{gather}
\begin{split}
& \begin{tikzpicture}[scale=0.35]
 \draw(0,-2) -- (0,0);
 \draw (0,0) arc [radius=0.5, start angle=180, end angle= 0];
 \draw(1,-1) -- (1,0);
 \draw(0.5,0.5) -- (0.5,1.5);
 \draw(-0.25,2.25) -- (-0.25,3);
 \draw [fill] (0.40,0.40) rectangle (0.60,0.60);
 \draw [fill] (-0.35,2.15) rectangle (-0.15,2.35);
 \draw (0.5,1.5) arc [radius=0.75, start angle=0, end angle= 180];
 \draw(-1,1.5) to [out=275,in=90] (-1.5,0);
 \draw[](1,-1) arc [radius=0.5, start angle=180, end angle= 270];
 \draw[ba](2,-1) arc [radius=0.5, start angle=360, end angle= 270];
 \draw(2,-1) -- (2,3);
 \draw[bam] (-0.25,3) arc [radius=1.125, start angle=180, end angle= 0];
 \draw [fill] (-1.6,-0.1) rectangle (-1.4,0.1);
 \draw (-2,-0.5) arc [radius=0.5, start angle=180, end angle= 0];
 \draw(-2,-0.5) -- (-2,-2);
 \draw(-1,-0.5) -- (-1,-2);
 \draw [fill] (-0.25, 3) circle[radius=0.15];
 \node at (4.35,1) {$=$};
 \draw [fill] (7.4, 2.4) rectangle (7.6,2.6);
 \draw (8.5,1.5) arc [radius=1, start angle=0, end angle= 180];
 \draw(6.5, 1.5) -- (6.5, -2);
 \draw(8.5, 1.5) -- (8.5, 1.25);
 \draw [fill] (8.4, 1.15) rectangle (8.6,1.35);
 \draw (9, 0.75) arc [radius=0.5, start angle=0, end angle= 180];
 \draw(8, 0.75) -- (8, 0);
 \draw [fill] (7.9, -0.1) rectangle (8.1,0.1);
 \draw (8.5, -0.5) arc [radius=0.5, start angle=0, end angle= 180];
 \draw(7.5, -0.5) -- (7.5, -2);
 \draw(8.5, -0.5) -- (8.5, -2);
 \draw(9, 0.75) -- (9, 0);
 \draw[] (9, 0) arc [radius=0.5, start angle=180, end angle= 360];
 \draw[ba](10,0) arc [radius=0.5, start angle=360, end angle= 270];
 \draw(10, 0) -- (10, 3);
 \draw[ibam] (10, 3) arc [radius=1.25, start angle=0, end angle= 180];
 \draw(7.5, 3) -- (7.5, 2.5);
 \draw [fill] (7.5, 3.2) circle [radius=0.15];
\node at (11.5,1) {$.$};
\end{tikzpicture}
\end{split}
\label{kappi-inv-pic}
\end{gather}
\end{prop}

\begin{proof}
We write the Killing form as in~\eqref{kappa_alt} and invoke the Jacobi identity in the form of~\eqref{jacobi_pic2} to
express the left hand side of~\eqref{kappi-inv-pic} as
\begin{gather}
\label{ki2}
\begin{split}
& \begin{tikzpicture}[scale=0.35]
 \draw(0,-2) -- (0,0);
 \draw (0,0) arc [radius=0.5, start angle=180, end angle= 0];
 \draw(1,-1) -- (1,0);
 \draw(0.5,0.5) -- (0.5,1.5);
 \draw(-0.25,2.25) -- (-0.25,3);
 \draw [fill] (0.40,0.40) rectangle (0.60,0.60);
 \draw [fill] (-0.35,2.15) rectangle (-0.15,2.35);
 \draw (0.5,1.5) arc [radius=0.75, start angle=0, end angle= 180];
 \draw(-1,1.5) to [out=275,in=90] (-1.5,0);
 \draw (1,-1) arc [radius=0.5, start angle=180, end angle= 270];
 \draw[ba] (2,-1) arc [radius=0.5, start angle=360, end angle= 270];
 \draw(2,-1) -- (2,3);
 \draw[bam] (-0.25,3) arc [radius=1.125, start angle=180, end angle= 0];
 \draw [fill] (-1.6,-0.1) rectangle (-1.4,0.1);
 \draw (-2,-0.5) arc [radius=0.5, start angle=180, end angle= 0];
 \draw(-2,-0.5) -- (-2,-2);
 \draw(-1,-0.5) -- (-1,-2);
 \draw [fill] (1, -0.6) circle[radius=0.15];
 \node at (4,1) {$=$};
 \draw [fill] (6.65,2.15) rectangle (6.85,2.35);
 \draw (7.5,1.5) arc [radius=0.75, start angle=0, end angle= 180];
 \draw(7.5,1.5) -- (7.5,1);
 \draw(6,1.5) -- (6,-2);
 \draw [fill] (7.4,0.9) rectangle (7.6,1.1);
 \draw (8.25,0.25) arc [radius=0.75, start angle=0, end angle= 180];
 \draw(6.75,0.25) -- (6.75,-2);
 \draw(8.25,0.25) -- (8.25,-0.25);
 \draw [fill] (8.15,-0.35) rectangle (8.35,-0.15);
 \draw (8.75,-0.75) arc [radius=0.5, start angle=0, end angle= 180];
 \draw(7.75,-0.75) -- (7.75,-2);
 \draw (8.75,-0.75) arc [radius=0.5, start angle=180, end angle= 270];
 \draw[ba] (9.75,-0.75) arc [radius=0.5, start angle=360, end angle= 270];
 \draw(9.75,-0.75) -- (9.75,3);
 \draw[ibam] (9.75,3) arc [radius=1.5, start angle=0, end angle= 180];
 \draw(6.75,3) -- (6.75,2.25);
 \draw [fill] (8.75, -0.8) circle[radius=0.15];
 \node at (11.5,1) {$-$};
 \draw [fill] (13.65,2.15) rectangle (13.85,2.35);
 \draw (14.5,1.5) arc [radius=0.75, start angle=0, end angle= 180];
 \draw(14.5,1.5) -- (14.5,1);
 \draw(13,1.5) to [out=275,in=90] (14,-2);
 \draw [fill] (14.4,0.9) rectangle (14.6,1.1);
 \draw (15.25,0.25) arc [radius=0.75, start angle=0, end angle= 180];
 \draw(13,-2) to [out=90,in=240 ] (13.5,-0.35);
 \draw(13.5,-0.35) to [out=60,in=275 ] (13.75,0.25);
 \draw(15.25,0.25) -- (15.25,-0.25);
 \draw [fill] (15.15,-0.35) rectangle (15.35,-0.15);
 \draw (15.75,-0.75) arc [radius=0.5, start angle=0, end angle= 180];
 \draw(14.75,-0.75) -- (14.75,-2);
 \draw (15.75,-0.75) arc [radius=0.5, start angle=180, end angle= 270];
 \draw[ba] (16.75,-0.75) arc [radius=0.5, start angle=360, end angle= 270];
 \draw(16.75,-0.75) -- (16.75,3);
 \draw[ibam] (16.75,3) arc [radius=1.5, start angle=0, end angle= 180];
 \draw(13.75,3) -- (13.75,2.25);
 \draw [fill] (15.75, -0.8) circle[radius=0.15];
 \node at (4+14+0.3,1) {$=$};
 \draw [fill] (6.65+14,2.15) rectangle (6.85+14,2.35);
 \draw (7.5+14,1.5) arc [radius=0.75, start angle=0, end angle= 180];
 \draw(7.5+14,1.5) -- (7.5+14,1);
 \draw(6+14,1.5) -- (6+14,-2);
 \draw [fill] (7.4+14,0.9) rectangle (7.6+14,1.1);
 \draw (8.25+14,0.25) arc [radius=0.75, start angle=0, end angle= 180];
 \draw(6.75+14,0.25) -- (6.75+14,-2);
 \draw(8.25+14,0.25) -- (8.25+14,-0.25);
 \draw [fill] (8.15+14,-0.35) rectangle (8.35+14,-0.15);
 \draw (8.75+14,-0.75) arc [radius=0.5, start angle=0, end angle= 180];
 \draw(7.75+14,-0.75) -- (7.75+14,-2);
 \draw (8.75+14,-0.75) arc [radius=0.5, start angle=180, end angle= 270];
 \draw[ba] (9.75+14,-0.75) arc [radius=0.5, start angle=360, end angle= 270];
 \draw(9.75+14,-0.75) -- (9.75+14,3);
 \draw[ibam] (9.75+14,3) arc [radius=1.5, start angle=0, end angle= 180];
 \draw(6.75+14,3) -- (6.75+14,2.25);
 \draw [fill] (8.75+14, -0.8) circle[radius=0.15];
 \node at (11.5+14,1) {$-$};
 \draw [fill] (13.65+14,2.15) rectangle (13.85+14,2.35);
 \draw (14.5+14,1.5) arc [radius=0.75, start angle=0, end angle= 180];
 \draw(14.5+14,1.5) -- (14.5+14,1);
 \draw(13+14,1.5) to (13+14,-2);
 \draw [fill] (14.4+14,0.9) rectangle (14.6+14,1.1);
 \draw (15.25+14,0.25) arc [radius=0.75, start angle=0, end angle= 180];
 \draw(15.65+14,-2) to [out=90,in=275 ] (13.75+14,0.25);
 \draw(15.25+14,0.25) -- (15.25+14,-0.25);
 \draw [fill] (15.15+14,-0.35) rectangle (15.35+14,-0.15);
 \draw (15.75+14,-0.75) arc [radius=0.5, start angle=0, end angle= 180];
 \draw(13.55+14,-2.2) to [out=90,in=275 ] (14.75+14,-0.75);
 \draw (15.75+14,-0.75) arc [radius=0.5, start angle=180, end angle= 270];
 \draw[ba] (16.75+14,-0.75) arc [radius=0.5, start angle=360, end angle= 270];
 \draw(16.75+14,-0.75) -- (16.75+14,3);
 \draw[ibam] (16.75+14,3) arc [radius=1.5, start angle=0, end angle= 180];
 \draw(13.75+14,3) -- (13.75+14,2.25);
 \draw [fill] (15.25+14, 0.3) circle [radius=0.15] ;
 \draw [fill] (14.5+14, -1.26) circle [radius=0.15] ;
 \node at (18+14,1){,};
\end{tikzpicture}
\end{split}
\end{gather}
where the second equality holds because of the cyclicity property~\eqref{eqcyclic}.
On the other hand, applying the Jacobi identity to the right hand side of~\eqref{kinv} yields
\begin{gather}
\label{222}
\begin{split}
& \begin{tikzpicture}[scale=0.35]
 \draw [fill] (-0.1, 2.4) rectangle (0.1,2.6);
 \draw (1,1.5) arc [radius=1, start angle=0, end angle= 180];
 \draw(-1, 1.5) -- (-1, -2);
 \draw(1, 1.5) -- (1, 1.25);
 \draw [fill] (0.9, 1.15) rectangle (1.1,1.35);
 \draw (1.5, 0.75) arc [radius=0.5, start angle=0, end angle= 180];
 \draw(0.5, 0.75) -- (0.5, 0);
 \draw [fill] (0.4, -0.1) rectangle (0.6,0.1);
 \draw (1, -0.5) arc [radius=0.5, start angle=0, end angle= 180];
 \draw(0, -0.5) -- (0, -2);
 \draw(1, -0.5) -- (1, -2);
 \draw(1.5, 0.75) -- (1.5, 0);
 \draw (1.5, 0) arc [radius=0.5, start angle=180, end angle= 270];
 \draw[->] (2.5,0) arc [radius=0.5, start angle=360, end angle= 270];
 \draw(2.5, 0) -- (2.5, 3);
 \draw[ibam] (2.5, 3) arc [radius=1.25, start angle=0, end angle= 180];
 \draw(0, 3) -- (0, 2.5);
 \draw [fill] (1.5, 0.25) circle [radius=0.15];
 \node at (4,1) {$=$};
 \draw [fill] (6.65,2.15) rectangle (6.85,2.35);
 \draw (7.5,1.5) arc [radius=0.75, start angle=0, end angle= 180];
 \draw(7.5,1.5) -- (7.5,1);
 \draw(6,1.5) -- (6,-2);
 \draw [fill] (7.4,0.9) rectangle (7.6,1.1);
 \draw (8.25,0.25) arc [radius=0.75, start angle=0, end angle= 180];
 \draw(6.75,0.25) -- (6.75,-2);
 \draw(8.25,0.25) -- (8.25,-0.25);
 \draw [fill] (8.15,-0.35) rectangle (8.35,-0.15);
 \draw (8.75,-0.75) arc [radius=0.5, start angle=0, end angle= 180];
 \draw(7.75,-0.75) -- (7.75,-2);
 \draw (8.75,-0.75) arc [radius=0.5, start angle=180, end angle= 270];
 \draw[->] (9.75,-0.75) arc [radius=0.5, start angle=360, end angle= 270];
 \draw(9.75,-0.75) -- (9.75,3);
 \draw[ibam] (9.75,3) arc [radius=1.5, start angle=0, end angle= 180];
 \draw(6.75,3) -- (6.75,2.25);
 \draw [fill] (8.75, -0.75) circle [radius=0.15];
 \node at (11.5,1) {$-$};
 \draw [fill] (13.65,2.15) rectangle (13.85,2.35);
 \draw (14.5,1.5) arc [radius=0.75, start angle=0, end angle= 180];
 \draw(14.5,1.5) -- (14.5,1);
 \draw(13,1.5) to (13,-2);
 \draw [fill] (14.4,0.9) rectangle (14.6,1.1);
 \draw (15.25,0.25) arc [radius=0.75, start angle=0, end angle= 180];
 \draw(15.65,-2) to [out=90,in=275 ] (13.75,0.25);
 \draw(15.25,0.25) -- (15.25,-0.25);
 \draw [fill] (15.15,-0.35) rectangle (15.35,-0.15);
 \draw (15.75,-0.75) arc [radius=0.5, start angle=0, end angle= 180];
 \draw(13.55,-2.2) to [out=90,in=275 ] (14.75,-0.75);
 \draw (15.75,-0.75) arc [radius=0.5, start angle=180, end angle= 270];
 \draw[->] (16.75,-0.75) arc [radius=0.5, start angle=360, end angle= 270];
 \draw(16.75,-0.75) -- (16.75,3);
 \draw[ibam] (16.75,3) arc [radius=1.5, start angle=0, end angle= 180];
 \draw(13.75,3) -- (13.75,2.25);
 \draw [fill] (15.75, -0.75) circle [radius=0.15];
\node at (18,1){.};
\end{tikzpicture}
\end{split}
\end{gather}
Analogously as in the proof of symmetry, the claim then follows noticing that the expressions on the right hand sides
of~\eqref{ki2} and~\eqref{222} are equal, owing to the functoriality of the twist.
\end{proof}

\subsection{Implications of non-degeneracy of the Killing form}
\label{sec3_4}

Non-degeneracy of the Killing form plays an important role for Lie algebras in categories of vector spaces.
It still does so for Lie algebras in generic additive symmetric ribbon categories.
We start by recalling the appropriate notion of non-degeneracy of a~pairing on an object~$U$, i.e.\ for a~morphism in
${\ensuremath{{\Hom}_{\ensuremath{\mathcal C}}}}(U\otimes U,{\bf1})$:

\begin{defi}
A~pairing $\varpi_U \in {\ensuremath{{\Hom}_{\ensuremath{\mathcal C}}}}(U\otimes U,{\bf1})$ in a~monoidal category
{\ensuremath{\mathcal C}}\ is called \emph{non-de\-ge\-nerate} if\/f there exists a~copairing $\varpi_U^- \in
{\ensuremath{{\Hom}_{\ensuremath{\mathcal C}}}}({\bf1},U\otimes U)$ that is side-inverse to $\varpi_U$, i.e.\ if\/f the
adjointness relations
\begin{gather}
\label{eq:sideinv}
(\varpi_U \otimes \id_U) \circ (\id_U \otimes \varpi_U^-) = \id_U = (\id_U \otimes \varpi_U) \circ (\varpi_U^- \otimes
\id_U)
\end{gather}
hold.
\end{defi}

We describe the equalities~\eqref{eq:sideinv} graphically as follows:
\begin{gather*}
\begin{tikzpicture} [scale=0.5]
 \draw(-2.25,-3.25) -- (-2.25,-1);
 \draw (-2.5,-1) rectangle (-1,-0.5+0.2);
 \node at (-1.75, -0.65) {$\scriptstyle{\varpi_{\scriptscriptstyle{U}}}$};
 \draw(-1.25,-2) -- (-1.25,-1);
 \draw[fill, fill opacity=0.1](-1.5,-2) rectangle (0,-2.5-0.2);
 \node at (-0.75, -2.3) {$\scriptstyle{\varpi^{-}_{\scriptscriptstyle{U}}}$};
 \draw(-0.25,-2) -- (-0.25,0.25);
 \node [below] at (-2.25, -3.5) {$\scriptstyle{U}$};
 \node [above] at (-0.25, 0.25) {$\scriptstyle{U}$};
 \node at (1, -1.5) {$=$};
 \draw(2,-3.25)--(2, 0.25);
 \node [below] at (2, -3.5) {$\scriptstyle{U}$};
 \node [above] at (2, 0.25) {$\scriptstyle{U}$};
 \node at (3, -1.5) {$=$};
 \draw[fill, fill opacity=0.1](4,-2) rectangle (5.5,-2.5-0.2);
 \node at (4.75, -2.3) {$\scriptstyle{\varpi^{-}_{\scriptscriptstyle{U}}}$};
 \draw(4.25,-2) -- (4.25,0.25);
 \draw (5,-1) rectangle (6.5,-0.5+0.2);
 \node at (5.75, -0.65) {$\scriptstyle{\varpi_{\scriptscriptstyle{U}}}$};
 \draw(5.25,-1) -- (5.25,-2);
 \draw(6.25,-1) -- (6.25,-3.25);
 \node [below] at (6.25, -3.5) {$\scriptstyle{U}$};
 \node [above] at (4.25, 0.25) {$\scriptstyle{U}$};
\end{tikzpicture}
\end{gather*}

If the pairing is symmetric, then each of the two equalities in~\eqref{eq:sideinv} implies the other.
If the monoidal category {\ensuremath{\mathcal C}}\ has a~right (say) duality, an equivalent def\/inition of
non-degeneracy of~$\varpi_U$ is that the morphism $\varpi_U^{\wedge}:= (\varpi_U\otimes \id_{U^\vee_{}}) \circ (\id_U
\otimes b_U)$ in ${\ensuremath{{\Hom}_{\ensuremath{\mathcal C}}}}(U,U^\vee)$ is an isomorphism, implying in particular
that $U^\vee$ is isomorphic, albeit not necessarily equal, to~$U$.

We have immediately the
\begin{lem}
\label{abkapp}
If a~Lie algebra has a~non-degenerate Killing form, then its only Abelian retract ideal is~$0$.
\end{lem}

\begin{proof}
Assume that $(K,e,r)$ is an Abelian retract ideal of a~Lie algebra~$L$.
Then, with $p = e \circ r$,
\begin{gather}
\kappa \circ (e \otimes \id_L)=\tracen3{p \circ \ell \circ (e \otimes \ell)} = \tracen3{\ell \circ (e\otimes\ell) \circ (\id_K \otimes \id_L \otimes p)}
\nonumber
\\
\phantom{\kappa \circ (e \otimes \id_L)}
= \tracen3{\ell \circ (e \otimes p) \circ (\id_K \otimes \ell) \circ (\id_K \otimes \id_L \otimes p)} = 0,
\label{333}
\end{gather}
where abelianness enters in the last step.
Graphically, with $e\,=$
\tikz[baseline=-0.75ex]{
   \draw[red](1,-0.5)--(1,0);
   \draw(1,0)--(1,0.4);
   \embed{1-0.075}{0}{0.15};}
\,and $r\,=$
\tikz[baseline=-0.75ex]{
   \draw(1,-0.5)--(1,0);
   \draw[red](1,0)--(1,0.4);
   \retract{1-0.075}{-0.125}{0.15};}\,,
equation~\eqref{333} reads
\begin{gather*}
  \begin{tikzpicture}[scale=0.3]
  \draw[red](-6-0.25,-2+0.5) -- (-6-0.25,-0.75+0.5);
  \draw(-6-0.25,-0.75+0.5) -- (-6-0.25,0.5+0.5);
  \node [left] at (-6-0.2-0.25,-0.85+0.5) {$\scriptstyle{e}$};
   \embed{-6-0.2-0.25}{-0.75+0.5}{0.4}
  \draw (-6.75-0.25,0.5+0.5) rectangle (-3.75-0.25,2+0.5);
   \node at (-5.25-0.25, 1.25+0.5) {$\scriptstyle{\kappa}$};
  \draw(-4.5-0.25,-2+0.5) -- (-4.5-0.25,0.5+0.5);
   \node [below] at (-6-0.25,-2+0.5) {$\scriptstyle{K}$};
  \node [below] at (-4.5-0.25,-2+0.5) {$\scriptstyle{L}$};
  \node at(-2.25-0.25,1){$\equiv$};
  \draw(0,-2) -- (0,0);
  \draw (0,0) arc [radius=0.5, start angle=180, end angle= 0];
  \draw(1,-1) -- (1,0);
  \draw(0.5,0.5) -- (0.5,1.5);
  \draw(-0.25,2.25) -- (-0.25,3);
  \draw [fill] (0.40,0.40) rectangle (0.60,0.60);
  \draw [fill] (-0.35,2.15) rectangle (-0.15,2.35);
  \draw [fill] (-0.25, 3) circle [radius=0.15] ;
  \draw (0.5,1.5) arc [radius=0.75, start angle=0, end angle= 180];
  \draw[red](-1.,0) -- (-1,-2);
   \draw(-1.,1.5) -- (-1,-0);
  \draw(1,-1) arc [radius=0.75, start angle=180, end angle= 270];
  \draw[ba]  (2.5,-1) arc [radius=0.75, start angle=360, end angle= 270];
  \draw(2.5,-1) -- (2.5,4);
  \draw(-0.25,2.5) -- (-0.25,4);
  \draw[ba] (-0.25,4) arc [radius=1.25+0.125, start angle=180, end angle= 90];
  \draw (1.25-0.25+0.125,4+1.25+0.125) arc [radius=1.25+0.125, start angle=90, end angle= 0];
  \embed{-1-0.2}{0}{0.4}
\node at (4,1){$=$} ;
  \draw(0+6.5,-2) -- (0+6.5,0);
  \draw (0+6.5,0) arc [radius=0.5, start angle=180, end angle= 0];
  \draw(1+6.5,-1) -- (1+6.5,0);
  \draw(0.5+6.5,0.5) -- (0.5+6.5,1.5);
  \draw(-0.25+6.5,2.25) -- (-0.25+6.5,3);
  \draw [fill] (0.40+6.5,0.40) rectangle (0.60+6.5,0.60);
  \draw [fill] (-0.35+6.5,2.15) rectangle (-0.15+6.5,2.35);
  \draw [fill] (-0.25+6.5, 2.75) circle [radius=0.15] ;
  \draw (0.5+6.5,1.5) arc [radius=0.75, start angle=0, end angle= 180];
   \draw[red](-1+6.5,0) -- (-1+6.5,-2);
   \draw(-1+6.5,1.5) -- (-1+6.5,-0);
    \draw(1+6.5,-1) arc [radius=0.75, start angle=180, end angle= 270];
  \draw[ba]  (2.5+6.5,-1) arc [radius=0.75, start angle=360, end angle= 270];
  \draw(2.5+6.5,-1) -- (2.5+6.5,4);
  \draw[red](-0.25+6.5,3.2) -- (-0.25+6.5,4);
  \retract{-0.25+6.5-0.2}{3}{0.4}
   \embed{-0.25+6.5-0.2}{4.25}{0.4}
  \draw[ba] (-0.25+6.5,4) arc [radius=1.25+0.125, start angle=180, end angle= 90];
  \draw (1.25-0.25+0.125+6.5,4+1.25+0.125) arc [radius=1.25+0.125, start angle=90, end angle= 0];
  \embed{-1-0.2+6.5}{0}{0.4}
\node at (10.5,1){$=$} ;
  \draw(0+13,-2) -- (0+13,0);
  \draw (0+13,0) arc [radius=0.5, start angle=180, end angle= 0];
  \draw[red](1+13,-1) -- (1+13,0);
  \retract{1+13-0.2}{-1.35}{0.4}
   \embed{1+13-0.2}{0}{0.4}
  \draw(0.5+13,0.5) -- (0.5+13,1.5);
  \draw(-0.25+13,2.25) -- (-0.25+13,3);
  \draw [fill] (0.40+13,0.40) rectangle (0.60+13,0.60);
  \draw [fill] (-0.35+13,2.15) rectangle (-0.15+13,2.35);
  \draw [fill] (-0.25+13, 3) circle [radius=0.15] ;
  \draw (0.5+13,1.5) arc [radius=0.75, start angle=0, end angle= 180];
   \draw[red](-1+13,0) -- (-1+13,-2);
   \draw(-1+13,1.5) -- (-1+13,-0);
    \draw(1+13,-1) arc [radius=0.75, start angle=180, end angle= 270];
  \draw[ba]  (2.5+13,-1) arc [radius=0.75, start angle=360, end angle= 270];
  \draw(2.5+13,-1) -- (2.5+13,4);
  \draw(-0.25+13,2.5) -- (-0.25+13,4);
  \draw[ba] (-0.25+13,4) arc [radius=1.25+0.125, start angle=180, end angle= 90];
  \draw (1.25-0.25+0.125+13,4+1.25+0.125) arc [radius=1.25+0.125, start angle=90, end angle= 0];
  \embed{-1-0.2+13}{0}{0.4}
\node at (17,1){$=$} ;
  \draw(0+19.5,-2) -- (0+19.5,-0.75);
  \draw (0+19.5,-0.75) arc [radius=0.5, start angle=180, end angle= 0];
  \draw(1+19.5,-1) -- (1+19.5,-0.75);
   \draw(0.5+19.5,-0.25) -- (0.5+19.5,0.25);
  \draw[red](0.5+19.5,0.25) -- (0.5+19.5,1.5);
  \retract{0.5-0.2+19.5}{0.25}{0.4}
  \embed{0.5-0.2+19.5}{1.5}{0.4}
 \draw(-0.25+19.5,2.25) -- (-0.25+19.5,3);
  \draw [fill] (0.40+19.5,0.40-0.75) rectangle (0.60+19.5,0.60-0.75);
  \draw [fill] (-0.35+19.5,2.15) rectangle (-0.15+19.5,2.35);
  \draw [fill] (-0.25+19.5, 3) circle [radius=0.15] ;
  \draw (0.5+19.5,1.5) arc [radius=0.75, start angle=0, end angle= 180];
   \draw[red](-1+19.5,0) -- (-1+19.5,-2);
   \draw(-1+19.5,1.5) -- (-1+19.5,-0);
    \draw(1+19.5,-1) arc [radius=0.75, start angle=180, end angle= 270];
  \draw[ba]  (2.5+19.5,-1) arc [radius=0.75, start angle=360, end angle= 270];
  \draw(2.5+19.5,-1) -- (2.5+19.5,4);
  \draw(-0.25+19.5,2.5) -- (-0.25+19.5,4);
  \draw[ba] (-0.25+19.5,4) arc [radius=1.25+0.125, start angle=180, end angle= 90];
  \draw (1.25-0.25+0.125+19.5,4+1.25+0.125) arc [radius=1.25+0.125, start angle=90, end angle= 0];
  \embed{-1-0.2+19.5}{0}{0.4}
\node at (24.5,1){$=\,\,0\,\,.$} ;
  \end{tikzpicture}
\end{gather*}
Using that $\kappa$ is non-degenerate, \eqref{333} implies immediately that
$0 = (\kappa\circ\mathrm{id}_L)\circ(e\otimes\kappa^{-})= e$.
\end{proof}

We are now in a~position to establish

\begin{prop}
\label{prop:main}
Let $(L,\ell)$ be a~Lie algebra in an idempotent complete symmetric ribbon category and let
${\ensuremath{{\End}_{\ensuremath{\mathcal C}}}}(L)$ have finitely many idempotents.
If the Killing form~$\kappa$ of~$L$ is non-degenerate, then~$L$ is a~finite direct sum of indecomposable Lie algebras.
\end{prop}

\begin{proof}
Assume that~$L$ is decomposable as a~Lie algebra, and let $(M,\e{}{},\r{}{})$ be a~non-trivial retract ideal of~$L$,
with corresponding idempotent $\er{}{} = \e{}{}\circ\r{}{}$.
The kernel of the idempotent $ \widehat{p}:=(\kappa\otimes\id_L)\circ \er{}{} \circ (\id_L\otimes \kappa^{-})$ exists
and is obtained by splitting the idempotent $p'=\id_L - \widehat{p} $; we denote this kernel retract by $(M', e', r')$,
so that
\begin{gather*}
\begin{tikzpicture} [scale=0.5]
\draw(-2.5,-2.5) -- (-2.5,-1);
\draw[blue](-2.5,-4) -- (-2.5,-2.5);
\embed{-2.5-0.15}{-2.5}{0.3}
\node [left] at (-2.5, -2.5) {$\scriptstyle{e'}$};
\node [below] at (-2.5, -4) {$\scriptstyle{M'}$};
\draw (-3,-1) rectangle (-1,-0.25);
\node at (-2, -0.65) {$\scriptstyle{\kappa}$};
\draw(-1.25,-1.25) -- (-1.25,-1);
\draw[red](-1.25,-1.25) -- (-1.25,-2.75);
\node [right] at (-1.4, -2) {$\scriptscriptstyle{M}$};
\embed{-1.25-0.15}{-1.25}{0.3}
\node [left] at (-1.25, -1.5) {$\scriptstyle{e}$};
\retract{-1.25-0.15}{-2.75}{0.3}
\node [left] at (-1.25, -2.6) {$\scriptstyle{r}$};
\draw(-1.25,-2.75) -- (-1.25,-3);
\draw[fill, fill opacity=0.2](-1.5,-3) rectangle (0.5,-3.75);
\node at (-0.5, -3.35) {$\scriptstyle{\kappa^{-}}$};
\draw(0,-3) -- (0,0.25);
\node [right] at (0, 0.25) {$\scriptstyle{L}$};
\node at (1.5, -1.5) {$=0.$};
\end{tikzpicture}
\end{gather*}
Since~$\kappa$ is non-degenerate, we have
\begin{gather}
\label{p1-1}
\kappa \circ  (e' \otimes e ) = 0.
\end{gather}
It follows in particular that
\begin{gather}
\label{p1-2}
0 = \kappa\circ (e'\otimes [p\circ \ell\circ (\id_{L} \otimes e)]) = \kappa \circ (e'\otimes [\ell\circ (\id_{L} \otimes e)])
= \kappa\circ ([\ell \circ (e'\otimes \id_{L})] \otimes e)
\end{gather}
or, graphically,
\begin{gather*}
\begin{tikzpicture} [scale=0.5]
\node at(-4.5,-0.85){$0=$};
\draw(-2.4,-3.25) -- (-2.4,-1);
\draw[blue](-2.4,-4.5) -- (-2.4,-3.25);
\embed{-2.4-0.15}{-3.25}{0.3}
\node [left] at (-2.2, -4.75) {$\scriptstyle{M'}$};
\draw (-3,-1) rectangle (-0.5,-0.25);
\node at (-1.75, -0.625) {$\scriptstyle{\kappa}$};
\draw(-1.1,-1.35) -- (-1.1,-1);
\draw[red](-1.1,-1.35) -- (-1.1,-2.5);
\embed{-1.1-0.15}{-1.25}{0.3}
\retract{-1.1-0.15}{-2.5}{0.3}
\draw(-1.1,-2.5) -- (-1.1,-2.9);
\draw(-1.6,-3.4) -- (-1.6,-4.5);
\draw (-1.6,-3.4) arc [radius=0.5, start angle=180, end angle= 0];
\draw(-0.6,-3.4) -- (-0.6,-3.65);
\draw [fill] (-1.05,-2.95) rectangle (-1.15,-2.85);
\draw[red](-0.6,-3.5) -- (-0.6,-4.5);
\embed{-0.6-0.15}{-3.35}{0.3}
\node [right] at (-0.85, -4.75) {$\scriptstyle{M}$};
\node [below] at (-1.6, -4.35) {$\scriptstyle{L}$};
\node at (0.5, -1) {$=$};
\draw(-2.4+4.5,-3.25) -- (-2.4+4.5,-1);
\draw[blue](-2.4+4.5,-4.5) -- (-2.4+4.5,-3.25);
\embed{-2.4-0.15+4.5}{-3.25}{0.3}
\draw (-3+4.5,-1) rectangle (-0.5+4.5,-0.25);
\node at (-1.75+4.5, -0.625) {$\scriptstyle{\kappa}$};
\draw(-1.1+4.5,-1.35) -- (-1.1+4.5,-1);
\draw(-1.1+4.5,-1.35) -- (-1.1+4.5,-2.5);
\draw(-1.1+4.5,-2.5) -- (-1.1+4.5,-2.9);
\draw(-1.6+4.5,-3.4) -- (-1.6+4.5,-4.5);
\draw (-1.6+4.5,-3.4) arc [radius=0.5, start angle=180, end angle= 0];
\draw(-0.6+4.5,-3.4) -- (-0.6+4.5,-3.65);
\draw [fill] (-1.05+4.5,-2.95) rectangle (-1.15+4.5,-2.85);
\draw[red](-0.6+4.5,-3.5) -- (-0.6+4.5,-4.5);
\embed{-0.6-0.15+4.5}{-3.35}{0.3}
\node at (5, -1) {$=$};
\draw(-2.3+4.5+5,-2.9) -- (-2.3+4.5+5,-1);
\draw[blue](-2.8+4.5+5,-4.5) -- (-2.8+4.5+5,-3.25);
\embed{-2.8-0.15+4.5+5}{-3.25}{0.3}
\draw (-3+4.5+5,-1) rectangle (-0.5+4.5+5,-0.25);
\node at (-1.75+4.5+5, -0.625) {$\scriptstyle{\kappa}$};
\draw(-1.1+4.5+5,-3.75) -- (-1.1+4.5+5,-1);
\draw(-1.8+4.5+5,-3.4) -- (-1.8+4.5+5,-4.5);
\draw (-1.8+4.5+5,-3.4) arc [radius=0.5, start angle=0, end angle= 180];
\draw [fill] (-2.35+4.5+5,-2.95) rectangle (-2.25+4.5+5,-2.85);
\draw[red](-1.1+4.5+5,-3.5) -- (-1.1+4.5+5,-4.5);
\embed{-1.1-0.15+4.5+5}{-3.35}{0.3}
\end{tikzpicture}
\end{gather*}
where in the second equality it is used that~$M$ is an ideal, and in the third that~$\kappa$ is invariant.
By the universal property of $M'$ as the kernel of $ \widehat{p}$,~\eqref{p1-2} implies
\begin{gather*}
\ell \circ (e' \otimes \id_{L}) = p' \circ \ell \circ (e' \otimes \id_{L})
\end{gather*}
showing that $(M', e', r')$ is in fact a~retract \emph{ideal} of~$L$.

Now take~$M$ to be a~\emph{primitive} retract ideal.
Then either $p\circ p' = p$ or $p \circ p' = 0$.
In case that $p \circ p' = p$, it follows from~\eqref{p1-1} that $\kappa \circ (e\otimes e) = 0$, which in turn implies
\begin{gather*}
\kappa \circ \big([\ell \circ (p \otimes p)]\otimes \id_L \big) = \kappa \circ \big(p \otimes [\ell\circ(p\otimes\id_L)] \big)
\\
\phantom{\kappa \circ \big([\ell \circ (p \otimes p)]\otimes \id_L \big)}{}
= \kappa \circ (p\otimes p) \circ \big(\id_L \otimes [\ell \circ(p\otimes\id_L)] \big) = 0
\end{gather*}
and therefore, when combined with Lemma~\ref{abkapp}, contradicts the non-degeneracy of~$\kappa$.
We thus conclude that $p \circ p' = 0$, and hence, again by non-degeneracy, that the Lie algebra~$L$ is the direct sum of
the retract ideals~$M$ and $M'$.
Furthermore, $\kappa\circ (e \otimes e)$ and $\kappa \circ (e' \otimes e')$ are just the Killing form of~$M$ and $M'$,
respectively.

The claim now follows by iteration, which terminates after a~f\/inite number of steps because
${\ensuremath{{\End}_{\ensuremath{\mathcal C}}}}(L)$ only has f\/initely many idempotents.
\end{proof}

\begin{rem}
A~large class of categories for which the assumptions of Proposition~\ref{prop:main} are satisf\/ied for any Lie algebra are Krull--Schmidt
categories, for which every object is a~f\/inite direct sum of indecomposables.
\end{rem}

\begin{rem}
For {\ensuremath{\mathcal C}}\ the category of f\/inite-dimensional vector spaces over a~f\/ield (of arbitrary
characteristic), the proof of Proposition~\ref{prop:main} reduces to the one given in~\cite{dieu}.
\end{rem}

As is well known, for Lie algebras in $\ensuremath{\mathscr V{\rm ect}}_\ensuremath{\Bbbk}$, semisimplicity does not imply
non-degeneracy of the Killing form unless \ensuremath{\Bbbk}\ is a~f\/ield of characteristic zero.
Still, the following weaker statement holds.

\begin{prop}
The Killing form of an indecomposable Lie algebra in a~symmetric ribbon category is either zero or non-de\-ge\-nerate.
\end{prop}

\begin{proof}
Let $(M, e, r)$ be a~maximal retract of the Lie algebra $(L,\ell)$ such that for the Killing form~$\kappa$ of~$L$ one
has
\begin{gather}
\label{kappML0}
\kappa \circ (\id_{L} \otimes e) = 0.
\end{gather}
Then we have
\begin{gather*}
0 = \kappa \circ (\ell \otimes e) = \kappa \circ (\id_{L} \otimes \ell) \circ (\id_{L} \otimes \id_{L} \otimes e),
\end{gather*}
which implies that~$M$ is a~retract ideal of~$L$.
Since~$L$ is indecomposable, this means that~$M$ either equals~$L$ or is zero.
In the former case,~\eqref{kappML0} says that~$\kappa$ is zero, while the latter case amounts to~$\kappa$ being
non-degenerate.
\end{proof}

\appendix

\section{Appendix}

\subsection{Categorical background}

For the sake of f\/ixing terminology and notation, we recall a~few pertinent concepts from category theory.
Besides in standard textbooks, details can e.g.\ be found on the~$n$Lab web site\footnote{\url{http://ncatlab.org}.}; for a~condensed exposition
of several of the relevant notions see~\cite{krau9}.

A \emph{subobject} $(U',e)$ of an object~$U$ is a~pair consisting of an object $U'$ and a~monomorphism $e\in
{\ensuremath{{\Hom}_{\ensuremath{\mathcal C}}}}(U',U)$ or, to be precise, an equivalence class of such pairs, with
equivalence def\/ined via factorization of morphisms with the same codomain (but it is common to use the term subobject
also for individual representatives.) A~\emph{simple} object is an object that does not have any non-trivial subobject.
A~\emph{retract} of an object~$U$ is a~subobject of~$U$ for which~$e$ has a~left-inverse, i.e.\ is a~triple $(U',e,r)$
such that $(U',e)$ is a~subobject of~$U$ and $r\in {\ensuremath{{\Hom}_{\ensuremath{\mathcal C}}}}(U',U)$ is an
epimorphism satisfying $r \circ e = \id_{U'}$.
A~retract $(U',e,r)$ for which the idempotent $p = e \circ r \in {\ensuremath{{\End}_{\ensuremath{\mathcal C}}}}(U)$
cannot be written as the sum of two non-zero idempotents is called \emph{primitive}.

The \emph{image} of a~morphism $f \in {\ensuremath{{\Hom}_{\ensuremath{\mathcal C}}}}(U,V)$ is a~pair $(V',h)$
consisting of an object $V'$ and a~monic $h \in {\ensuremath{{\Hom}_{\ensuremath{\mathcal C}}}}(V',V)$ such that there
exists a~morphism $g \in {\ensuremath{{\Hom}_{\ensuremath{\mathcal C}}}}(U,V')$ satifying $h \circ g = f$ and such that
$(V',h)$ is universal with this property, i.e.\ for any object~$W$ with a~monic $k \in
{\ensuremath{{\Hom}_{\ensuremath{\mathcal C}}}}(W,V)$ and morphism $j \in {\ensuremath{{\Hom}_{\ensuremath{\mathcal
C}}}}(U,W)$ satifying $k \circ j = f$ there exists a~unique morphism $l \in {\ensuremath{{\Hom}_{\ensuremath{\mathcal
C}}}}(V',W)$ such that $k \circ l = h$.
In short, in case it exists, the image is the smallest subobject $(V',h)$ of~$V$ through which the morphism~$f$ factors.

A \emph{preadditive} category is a~category {\ensuremath{\mathcal C}}\ which is enriched over the monoidal category of
Abelian groups, i.e.\ for which the morphism set ${\ensuremath{{\Hom}_{\ensuremath{\mathcal C}}}}(U,V)$ for any pair of
objects $U,V\in {\ensuremath{\mathcal C}}$ is an Abelian group and composition of morphisms is bilinear.
An \emph{additive} category is a~preadditive category admitting all f\/initary biproducts, or direct sums.
A~direct sum is unique up to isomorphism, and the empty biproduct is a~zero object 0 satisfying
${\ensuremath{{\Hom}_{\ensuremath{\mathcal C}}}}(U,0) = 0 = {\ensuremath{{\Hom}_{\ensuremath{\mathcal C}}}}(0,U)$ for
any object $U\in {\ensuremath{\mathcal C}}$.
An object~$U$ of an additive category is called \emph{indecomposable} if\/f it is non-zero and has only trivial direct sum
decompositions, i.e.\ $U = U_1 \oplus  U_2$ implies that $U_1 = 0$ or $U_2 = 0$.
Each summand $U_i$ in a~direct sum $U = U_1 \oplus U_2 \oplus \dots \oplus U_n$ forms a~retract $(U_i,e_i,r_i)$, i.e.\
there are monics $e_i \in {\ensuremath{{\Hom}_{\ensuremath{\mathcal C}}}}(U_i,U)$ and epis $r_i \in
{\ensuremath{{\Hom}_{\ensuremath{\mathcal C}}}}(U,U_i)$ such that
\begin{gather*}
r_i \circ e_j = \delta_{i,j} \id_{U_j}^{}
\qquad
\text{for~all}
\quad
i, j = 1,2,\dots,n
\qquad
\text{and}
\qquad
\sum\limits_{i=1}^n e_i \circ r_i = \id_U.
\end{gather*}
An indecomposable object thus only has trivial retracts, the zero object and itself.
A~primitive retract is an indecomposable direct summand.

An \emph{idempotent}~$p$ is an endomorphism satisfying $p^2 = p$; an idempotent $p \in
{\ensuremath{{\End}_{\ensuremath{\mathcal C}}}}(U)$ is called \emph{split} if\/f there exists a~retract $(V,e,r)$ of~$U$
such that $e \circ r = p$.
A~category {\ensuremath{\mathcal C}}\ is said to be \emph{idempotent complete} if\/f every idempotent in
{\ensuremath{\mathcal C}}\ is split.
In an idempotent complete category, any retract is a~direct sum of primitive retracts.
A~\emph{Krull--Schmidt category} is an additive category for which every object decomposes into a~f\/inite direct sum of
objects whose endomorphism sets are local rings.
The objects with local endomorphism rings are then precisely the indecomposable objects.
Krull--Schmidt categories are idempotent complete.
The class of Krull--Schmidt categories contains in particular every Abelian category in which each object has f\/inite
length, and thus e.g.\ all f\/inite tensor categories in the sense of~\cite{etos}, as well as~\cite{happ3} every
triangulated category whose morphism sets are f\/inite-dimensional vector spaces and for which the endomorphism ring of
any indecomposable object is local.

The data of a~\emph{monoidal} category $({\ensuremath{\mathcal C}},\otimes,{\bf1},a,l,r)$ consist of a~category
{\ensuremath{\mathcal C}}, the tensor product functor $\otimes\colon {\ensuremath{\mathcal C}} {\times}
{\ensuremath{\mathcal C}} {\to} {\ensuremath{\mathcal C}}$, the monoidal unit ${\bf1} \in {\ensuremath{\mathcal C}}$,
the associativity constraint $a = (a_{U,V,W})$, and left and right unit constraints $l = (l_U)$ and $r = (r_U)$ (see
e.g.~\cite{muge18} for a~review).
If {\ensuremath{\mathcal C}}\ is additive, then the tensor product is required to be additive in each argument.
Every monoidal category is equivalent to a~\emph{strict} one, i.e.\ one for which the associativity and unit constraints
are identities.
We denote by $U^{\otimes n}$ the tensor product of~$n$ copies of~$U$ (which can be def\/ined unambiguously even if
{\ensuremath{\mathcal C}}\ isn't strict).
A~\emph{braided} monoidal category has in addition a~commutativity constraint, called the braiding and denoted by $c =
(c_{U,V})$, i.e.\ a~natural family of isomorphisms $c_{U,V} \in {\ensuremath{{\Hom}_{\ensuremath{\mathcal C}}}}(U\otimes
V,V\otimes U)$ satisfying two coherence conditions known as hexagon equations.
A~\emph{symmetric} monoidal category is a~braided monoidal category with symmetric braiding, i.e.\ for which $c_{V,U}
\circ c_{U,V} = \id_{U\otimes V}$.
A~\emph{rigid} (or autonomous) monoidal category has right and left dualities, with natural families $d = (d_U)$ and
$\tilde d = (\tilde d_U)$ of evaluation and $b = (b_U)$ and $\tilde b = (\tilde b_U)$ of coevaluation morphisms,
respectively.
The objects that are right and left dual to~$U$ are denoted by $U^\vee$ and ${}^{\vee }U$, respectively, so that $d_U
\in {\ensuremath{{\Hom}_{\ensuremath{\mathcal C}}}} (U^{\vee} \otimes U, {\bf1})$, $b_U \in
{\ensuremath{{\Hom}_{\ensuremath{\mathcal C}}}}({\bf1}, U \otimes U^{\vee})$ and $\tilde{d}_U \in
{\ensuremath{{\Hom}_{\ensuremath{\mathcal C}}}}(U \otimes {}^{\vee}U,{\bf1})$, $\tilde{b}_U \in
{\ensuremath{{\Hom}_{\ensuremath{\mathcal C}}}}({\bf1},{}^{\vee}U \otimes U)$.
(With this convention, the functor $U^\vee{\otimes} -$ is left adjoint to $U \otimes -$, while ${}^{\vee}U \otimes -$ is
right adjoint to $U \otimes -$.) Dualities are also def\/ined naturally on morphisms, namely for $f\in
{\ensuremath{{\Hom}_{\ensuremath{\mathcal C}}}}(U,V)$~by
\begin{gather*}
\begin{tikzpicture}[scale=0.45]
\node at (-11.95,0){$f^\vee := (d_V \otimes \id_{U^\vee_{}}) \circ (\id_{V^\vee_{}} \otimes f \otimes \id_{U^\vee_{}})
\circ (\id_{V^\vee_{}} \otimes b_U)=$};
\draw[bam] (0.5,1) arc [radius=0.5, start angle=0, end angle= 180];
\draw (0.5,0.35)--(0.5,1);
\draw (-0.5,1)--(-0.5,-1.75);
\draw (0,-0.35) rectangle (1,0.35);
\node at (0.5,0) {$\scriptstyle{f}$};
\draw (0.5,-0.35)--(0.5,-1);
\draw[ibam] (0.5,-1) arc [radius=0.5, start angle=180, end angle= 360];
\draw (1.5,1.75)--(1.5,-1);
\node[left] at (-0.35,-1.65) {$\scriptscriptstyle{V^{\vee}}$};
\node[right] at (1.5,1.65) {$\scriptscriptstyle{U^{\vee}}$};
\node at (3,0) {$$};
\end{tikzpicture}
\end{gather*}
and
\begin{gather*}
\begin{tikzpicture}[scale=0.45]
\node at (-11.65,0){${}^{\vee\!\!}f:= (\id_{{}_{}^{\vee\!}U} \otimes \tilde d_V) \circ (\id_{{}_{}^{\vee\!}U}
\otimes f \otimes \id_{{}_{}^{\vee\!}V}) \circ (\tilde b_U {\otimes} \id_{{}_{}^{\vee\!}V})=$};
\draw[bam] (0.5,1) arc [radius=0.5, start angle=180, end angle= 0];
\draw (0.5,0.35)--(0.5,1);
\draw (1.5,1)--(1.5,-1.75);
\draw (0,-0.35) rectangle (1,0.35);
\node at (0.5,0) {$\scriptstyle{f}$};
\draw (0.5,-0.35)--(0.5,-1);
\draw[ibam] (0.5,-1) arc [radius=0.5, start angle=360, end angle= 180];
\draw (-0.5,-1.01)--(-0.5,1.75);
\node[left] at (-0.35,1.65) {$\scriptscriptstyle{^{\vee}V}$};
\node[right] at (1.5,-1.65) {$\scriptscriptstyle{^{\vee}U}$};
\node at (3,0) {$,$};
\end{tikzpicture}
\end{gather*}
respectively, whereby they furnish (contravariant) endofunctors of {\ensuremath{\mathcal C}}.
For a~monoidal category with a~right (say) duality, $(U,e,r)$ is a~retract of~$V$ if\/f $(U^\vee,r^\vee,e^\vee)$ is
a~retract of $V^\vee$.

A \emph{sovereign} structure on a~rigid monoidal category is a~(choice of) monoidal natural isomorphism between the
right and left duality functors, i.e.\ a~natural family of isomorphisms $\sigma_U \in
{\ensuremath{{\Hom}_{\ensuremath{\mathcal C}}}}(U^\vee_{},{}_{}^{\vee }U)$ such that
\begin{gather}
{}^{\vee }f \circ \sigma_V = \sigma_U \circ f^\vee
\label{fvvf}
\end{gather}
for all $f\in {\ensuremath{{\Hom}_{\ensuremath{\mathcal C}}}}(U,V)$.
A~sovereign structure is equivalent~\cite[Proposition~2.11]{yett2} to a~pivotal (or balanced) structure, i.e.\ to a~monoidal natural isomorphism
between the (left or right) double dual functor and the identity functor.
If the natural isomorphism def\/ining a~sovereign structure is the identity, the category is called \emph{strictly sovereign}.
In a~sovereign category, an endomorphism $f\in {\ensuremath{{\End}_{\ensuremath{\mathcal C}}}}(U)$ has a~right trace
$\tracen{}{f} \in {\ensuremath{{\End}_{\ensuremath{\mathcal C}}}}({\bf1})$, given (for strictly sovereign
{\ensuremath{\mathcal C}}) by $ \tracen{}{f} = \tilde{d}_U \circ (f\otimes\id_{U_{}^{\vee}})\circ b_{U}$, as well as an
analogous left trace, and thus in particular every object~$U$ has a~right (and analogously, left) dimension
${\ensuremath{{\dim}_{\ensuremath{\mathcal C}}}}(U)  \equiv  \tracen{}{\id_U} \in
{\ensuremath{{\End}_{\ensuremath{\mathcal C}}}}({\bf1})$.
When the left and right traces coincide, the category is called \emph{spherical}.

A \emph{twist} (or balancing) for a~braided monoidal category with a~(right, say) duality is a~na\-tu\-ral family $\theta =
(\theta_U)$ of isomorphisms $\theta_U \in {\ensuremath{{\End}_{\ensuremath{\mathcal C}}}}(U)$ satisfying $\theta_{\bf1}
= \id_{\bf1}$, $\theta_{U\otimes V} = (\theta_U \otimes\theta_V) \circ c_{V,U} \circ c_{U,V}$ and $\theta_{U^\vee_{}} =
{(\theta_U)}^\vee$.
A~braided rigid monoidal category equipped with a~compatible twist is called a~\emph{ribbon} (or tortile) category.
A~ribbon category is sovereign, with sovereign structure
\begin{gather}
\label{defsigma}
\sigma_U = (\id_{{}_{}^{\vee }U} \otimes d_U) \circ \big(c_{U,U}^{-1} \otimes \theta_U^{}\big) \circ (\id_{U^\vee_{}} \otimes
\tilde b_U) \in {\ensuremath{{\Hom}_{\ensuremath{\mathcal C}}}}(U^\vee_{}, {}_{}^{\vee }U),
\end{gather}
as well as spherical.
In a~symmetric ribbon category the twist $\theta_U$ squares to $\id_U$ (see e.g.~\cite[Remark 4.32]{seli3}).

\subsection{Algebras and coalgebras}

We also collect a~few elementary concepts concerning algebras in monoidal categories.
An \emph{algebra} (or monoid, or monoid object) in a~monoidal category ${\ensuremath{\mathcal C}} =
({\ensuremath{\mathcal C}},\otimes,{\bf1},a,l,r)$ is just a~pair $(A,m)$ consisting of an object $A \in
{\ensuremath{\mathcal C}}$ and a~morphism $m \in {\ensuremath{{\Hom}_{\ensuremath{\mathcal C}}}}(A\otimes A,A)$;~$m$ is
called a~product or a~multiplication, and by the usual abuse of terminology also the object~$A$ is referred to as an algebra.

This concept is rather empty unless an interesting property of the product and/or additional structure on~$A$ is
imposed, which in turn may require also the category {\ensuremath{\mathcal C}}\ to have additional structure beyond
being monoidal.
Monoidal categories without any further structure furnish the proper setting for \emph{associative} algebras, for which
the associativity property
\begin{gather}
m \circ (m \otimes \id_{A}) = m \circ (\id_{A} \otimes m) \circ a^{}_{A,A,A}
\label{assoc}
\end{gather}
is imposed, and for \emph{unital} algebras, for which there exists a~morphism $\eta \in
{\ensuremath{{\Hom}_{\ensuremath{\mathcal C}}}}({\bf1},A)$ satisfying
\begin{gather}
m \circ (\eta \otimes\ \id_A) \circ l_A^{-1} = \id_A = m \circ (\id_A \otimes \eta) \circ r_{ A}^{-1}.
\label{unital}
\end{gather}
Likewise, in any monoidal category there is the notion of a~coassociative and co-unital coalgebra
$(C,\Delta,\varepsilon)$, with morphisms $\Delta \in {\ensuremath{{\Hom}_{\ensuremath{\mathcal C}}}}(C,C\otimes C)$ and
$\varepsilon \in {\ensuremath{{\Hom}_{\ensuremath{\mathcal C}}}}(C, {\bf1})$ that satisfy relations obtained by
reversing the arrows in~\eqref{assoc} and~\eqref{unital}, i.e.\ $ (\Delta \otimes \id_{C}) \circ \Delta \circ a^{}_{C,C,C}
= (\id_{C} \otimes \Delta) \circ \Delta$ and $ l_C \circ (\varepsilon \otimes \id_C) \circ \Delta = \id_C = r_C \circ (\id_C
\otimes \varepsilon) \circ \Delta$.
In the sequel we will simplify such equalities by taking the monoidal category {\ensuremath{\mathcal C}}\ to be strict.

Still working merely in a~(strict) monoidal category {\ensuremath{\mathcal C}}\ there is an interesting way to combine
algebras and coalgebras:

\begin{defi}\label{def:Frob}\qquad
\begin{enumerate}\itemsep=0pt
\item[(i)] A~\emph{Frobenius algebra}~$A$ in {\ensuremath{\mathcal C}}\ is a~quintuple $(A,m,\eta,\Delta,\varepsilon)$
such that $(A,m,\eta)$ is a~unital associative algebra, $(A,\Delta,\varepsilon)$ is co-unital coassociative coalgebra,
and such that the coproduct is a~morphism of~$A$-bimodules, i.e.~
\begin{gather*}
(\id_{A}\otimes m) \circ (\Delta\otimes\id_{A}) = \Delta \circ m = (m \otimes\id_{A}) \circ (\id_{A} \otimes \Delta).
\end{gather*}
\item[(ii)] A~\emph{strongly separable}\footnote{In the literature also variants of this def\/inition are in use.}
Frobenius algebra is a~Frobenius algebra for which the product is a~left inverse of the coproduct, i.e.~$m \circ \Delta =
\id_{A}$.
\end{enumerate}
\end{defi}

\begin{Example}\quad
\begin{enumerate}\itemsep=0pt
\item[(i)] Coalgebras, and in particular Frobenius algebras, in symmetric monoidal categories which admit the notion of an
adjoint (or dagger~\cite{seli}) of a~morphism arise naturally in a~categorical formalization of fundamental issues of
quantum mechanics, like complementarity of observables or the no-cloning theorem; see e.g.~\cite{coDu2,copV,kissi3}.
\item[(ii)] To a~pair of functors $F\colon {\ensuremath{\mathcal C}} {\to} \mathcal D$ and $G\colon \mathcal D {\to}
{\ensuremath{\mathcal C}}$, such that~$G$ is right adjoint to~$F$, there is naturally associated a~monad in
{\ensuremath{\mathcal C}}, i.e.\ an algebra in the category of endofunctors of~{\ensuremath{\mathcal C}}.
If~$G$ is also left adjoint to~$F$, then this monad has a~natural Frobenius structure~\cite[Proposition~1.5]{stre8}.
\end{enumerate}
\end{Example}

If the monoidal category {\ensuremath{\mathcal C}}\ is in addition braided, one can also def\/ine (co-)commutativity as
well as bi- and Hopf algebras.
An algebra $(A,m)$ in a~braided (strict) monoidal category $({\ensuremath{\mathcal C}},\otimes,{\bf1},c)$ is said to be
\emph{commutative} if\/f $m\circ c_{A,A} = m$; dually, a~coalgebra $(C,\Delta)$ in {\ensuremath{\mathcal C}}\ is
\emph{co\-com\-mutative} if\/f $c_{A,A} \circ \Delta = \Delta$.
A~\emph{bialgebra}~$B$ is both an (associative unital) algebra and a~(co\-as\-sociative co-unital) coalgebra, with the
coproduct and counit being algebra morphisms, i.e.~$\Delta\circ m= (m\otimes m)\circ (\id_{B}\otimes c_{B,B}\otimes
\id_{B})\circ (\Delta\otimes\Delta)$ and $\varepsilon\circ m= \varepsilon\otimes\varepsilon$.
A~\emph{Hopf} algebra~$H$ is a~bialgebra with an antipode $s \in {\ensuremath{{\End}_{\ensuremath{\mathcal C}}}}(H)$
that is a~two-sided inverse of $\id_H$ for the convolution product.

\begin{Example}\label{exa:4}\quad
\begin{enumerate}
\itemsep=0pt
\item[(i)] If {\ensuremath{\mathcal C}}\ is braided monoidal and countable direct sums are def\/ined in {\ensuremath{\mathcal
C}},\footnote{If {\ensuremath{\mathcal C}}\ is Abelian, then the condition that countable direct sums exist is
dispensable: instead of {\ensuremath{\mathcal C}}\ one may then work with the category of ind-objects in
{\ensuremath{\mathcal C}}, into which {\ensuremath{\mathcal C}}\ is fully embedded~\cite[Section~3.1]{hai2}.} then for
any object $U \in {\ensuremath{\mathcal C}}$ the object $T(U):= \bigoplus_{n=0}^\infty U^{\otimes n}$ carries a~natural
structure of a~unital associative algebra, with the product def\/ined via the tensor product of the retracts $U^{\otimes
n}$ and the unit morphism given by the retract embedding of $U^{\otimes 0} = {\bf1}$, see e.g.~\cite{DeM, schau6};
$T(U)$ is called the \emph{tensor algebra} of~$U$.
In fact~\cite[Corollary~2.4]{schau6}, $T(U)$ has in fact a~natural \emph{bi}algebra structure, and by extending the endomorphism $-\id_U$ of~$U$
to an anti-algebra morphism of $T(U)$ def\/ines an antipode, so that $T(U)$ is even a~Hopf algebra.
$T(U)$ can also be def\/ined by a~universal property, it is an associative algebra together with a~morphism $i_U \in
{\ensuremath{{\Hom}_{\ensuremath{\mathcal C}}}}(U,T(U))$ such that for any algebra morphism~$f$ from~$U$ to some
algebra~$A$ in~{\ensuremath{\mathcal C}}\ there exists a~unique morphism~$g$ from $T(U)$ to~$A$ such that $f = g \circ
i_U$.
\item[(ii)] If in addition {\ensuremath{\mathcal C}}\ is idempotent complete, is $\ensuremath{\Bbbk}$-linear with
$\ensuremath{\Bbbk}$ a~f\/ield of characteristic zero, and has symmetric braiding, then the direct sum $S(U):=
\bigoplus_{n=0}^\infty S_n(U)$ of all symmetric tensor powers of~$U$ carries a~natural structure of a~commutative unital
associative algebra.
\item[(iii)] Hopf algebras appear naturally in constructions within categories of three-dimensional cobordisms.
Specif\/ically, the manifold obtained from a~2-torus by excising an open disk is a~Hopf algebra 1-morphism in the
bicategory of three-dimensional cobordisms with corners~\cite{crye2}.
\item[(iv)] For a~braided f\/inite tensor category {\ensuremath{\mathcal C}}\ the coend $\int^{U\in{\ensuremath{\mathcal
C}}} U^\vee{\otimes} U$ exists and can be used to construct quantum invariants of three-manifolds and representations of
mapping class groups~\cite{lyub8}.
It also plays a~role (for {\ensuremath{\mathcal C}}\ semisimple, i.e.\ a~fusion category) in establishing the
relationship between Reshetikhin--Turaev and Turaev--Viro invariants~\cite{tuVi}.
For these results it is crucial that this coend carries a~natural structure of a~Hopf algebra (with some additional
properties) in {\ensuremath{\mathcal C}}.
\end{enumerate}
\end{Example}

\begin{rem}
A~Hopf algebra~$H$ in an additive rigid braided category which has invertible antipode and which possesses a~left
integral $\Lambda \in {\ensuremath{{\Hom}_{\ensuremath{\mathcal C}}}}({\bf1},H)$ and right cointegral $\lambda \in
{\ensuremath{{\Hom}_{\ensuremath{\mathcal C}}}}(H,{\bf1})$ such that $\lambda \circ \Lambda \in
{\ensuremath{{\End}_{\ensuremath{\mathcal C}}}}({\bf1})$ is invertible, naturally also has the structure of a~Frobenius
algebra with the same algebra structure, and with the Frobenius counit given by~$\lambda$~\cite{fuSc17,kaSt2, pare7}.
\end{rem}

Frobenius algebras~$F$ can, just like in the classical case, also be def\/ined with the help of a~non-degenerate invariant
pairing $\varpi \in {\ensuremath{{\Hom}_{\ensuremath{\mathcal C}}}}(F\otimes F,{\bf1})$; in terms of the description
chosen above, $\varpi = \varepsilon \circ m$, see e.g.~\cite{fuSt, yama12}.
A~\emph{symmetric} Frobenius algebra is a~Frobenius algebra~$F$ for which the pairing~$\varpi$ is symmetric.

\begin{Example}
In two-dimensional rational conformal f\/ield theory (RCFT), full local RCFTs that come from one and the same chiral RCFT
are in bijection with Morita classes of strongly separable symmetric Frobenius algebras of non-zero dimension in the
modular tensor category that underlies the chiral RCFT~\cite{fuRs4}.
\end{Example}

\subsection*{Acknowledgments}
We are most grateful to the referees for their valuable comments on an earlier version of this note, and in particular
for suggesting improvements of the proof of Proposition~\ref{prop:main}.
We also thank Christoph Schweigert for discussions and Scott Morrison for bringing Example~\ref{ex:M..}(iv) to our attention.
JF is supported by VR under project no.~621-2013-4207.
JF thanks the Erwin-Schr\"odinger-Institute (ESI)
for the hospitality during the programs ``Modern Trends in TQFT'' and
``Topological Phases of Quantum Matter'' while part of this work was pursued.

\pdfbookmark[1]{References}{ref}
\LastPageEnding

\end{document}